\documentclass[12pt]{article}
\usepackage{latexsym}
\usepackage{amsfonts}
\usepackage{amssymb}
\newtheorem{Th}{Theorem}[section]
\newtheorem{Def}{Definition}[section]
\newtheorem{Prop}{Proposition}[section]
\newtheorem{Le}{Lemma}[section]
\newtheorem{Cor}{Corollary}[section]
\newtheorem{Rem}{Remark}[section]

\newcommand\R{{\mathbb{R}}}
\newcommand\C{{\mathbb{C}}}

\newcommand\Fc{\mathcal{F}}
\newcommand\Hc{\mathcal{H}}
\newcommand\Lc{\mathcal{L}}
\newcommand\Mc{\mathcal{M}}
\newcommand\Rc{\mathcal{R}}
\newcommand\Sc{\mathcal{S}}
\newcommand\Tc{\mathcal{T}}
\newcommand\Wc{\mathcal{W}}
\newcommand\kn{\bigcirc\hspace{-0,45cm}\wedge\ } 
\newcommand\Zi{Z_{i}}
\newcommand\Zj{Z_{j}}
\newcommand\ftn{\phi^{*}TN}
\newcommand\Ftn{\Phi^{*}TN}
\newcommand\wx[1]{{\wedge}_{\xi}^{#1}(M)}
\newcommand\wh[1]{{\wedge}_{H}^{#1}(M)}
\newcommand\who{{\wedge}_{H_{0}}^{2}(M)}
\newcommand\whpm[1]{{\wedge}_{H}^{2,#1}(M)}
\newcommand\whtpm[1]{{\wedge}_{H_{#1}}^{2}(M)}
\newcommand\sh[1]{{S}_{H}^{#1}(M)}
\newcommand\ox[1]{{\Omega}_{\xi}^{#1}(M)}
\newcommand\oh[1]{{\Omega}_{H}^{#1}(M)}
\newcommand\ohe[1]{{\Omega}_{H}^{#1}(M;E)}

\newcommand\ep[1]{\epsilon_{#1}}
\newcommand\qhc{\stackrel{\circ}{Q_{H}}}
\newcommand\qhh{\widehat{{Q}_{H}}}
\newcommand\qh[4]{{Q}_{H}(\ep{#1},\ep{#2},\ep{#3},\ep{#4})}

\newcommand\sgh{{\sigma}_{H}}

\newcommand\dfh{d{\phi}_{H}}

\newcommand\dfxh{{(d\phi(\xi))}_{H^{'}}}
\newcommand\dfhh{d{\phi}_{H,H^{'}}}
\newcommand\nbsh{{\nabla}^{\Sc}\sgh}
\newcommand\nbfh{{\nabla}^{\Sc}\dfh}
\newcommand\nbfho{{\nabla}^{\Sc}d\phi_{H_{0}}}
\newcommand\delfh{\delta^{\nabla^{'}}_{H}\dfh}
\newcommand\dhsh{{d}_{H}^{\nabla^{E}}\sgh}
\newcommand\dhfh{{d}_{H}^{\nabla^{'}}\dfh}
\newcommand\nbish[2]{({\nabla}^{\Sc}\sgh)(\ep{#1},\ep{#2})}
\newcommand\dhish[2]{({d}_{H}^{\nabla^{E}}\sgh)(\ep{#1},\ep{#2})}
\newcommand\lx{{\Lc}_{\xi}^{\nabla^{E}}}
\newcommand\nab[1]{{\nabla}_{\ep{#1}}}
\newcommand\nabx{\nabla_{\xi}}
\newcommand\delhj{\delta^{\nabla^{'}}_{H,J}}
\newcommand\nn[1]{|{#1}|^{2}}

\newcommand\rw{R^{\Wc}_{H}}
\newcommand\rhp{R_{H^{'}}}
\newcommand\rwp{R^{\Wc}_{H^{'}}}
\newcommand\rwpf{{(\phi^{*}\rwp)}_{H}}

\newcommand\rwpfh[1]{\widehat{{(\phi^{*}\rwp)}_{H}^{#1}}}
\newcommand\rwpfhc{\widehat{{(\phi^{*}\rwp)}_{H}^{\C}}}
\newcommand\Cm{C^{\Mc}_{H}}
\newcommand\rwo{R^{\Wc}_{H_{0}}}
\newcommand\riw{Ric^{\Wc}_{H}}
\newcommand\riwo{Ric^{\Wc}_{H_{0}}}
\newcommand\row{{\rho}^{\Wc}_{H}}
\newcommand\rowo{{\rho}^{\Wc}_{H_{0}}}
\newcommand\sw{s^{\Wc}}
\newcommand\rf[2]{tr^{#1,#2}_{H}(\rwpfhc)}

\begin{document}

\title{\LARGE\bf Mok-Siu-Yeung type formulas on contact locally sub-symmetric spaces}	   
\author{Robert PETIT*} 
\date{}
\maketitle
\vspace{2cm}
\begin{abstract}
We derive Mok-Siu-Yeung type formulas for horizontal maps from compact contact locally sub-symmetric spaces into strictly pseudoconvex $CR$ manifolds and we obtain some rigidity theorems for the horizontal pseudoharmonic maps. 
\end{abstract}

\vspace{3cm}\noindent
{\bf Author Keywords:} Contact locally sub-symmetric spaces, Tanaka-Webster connection, pseudoharmonic maps, CR maps, Rumin complex, Mok-Siu-Yeung type formulas, rigidity results.\\

\noindent
{\bf Mathematical subject codes:} 32V05, 32V10, 53C17, 53C21, 53C24, 53C25, 53C30, 53C35, 53D10, 58E20\\

\vspace{1cm}\noindent
*Laboratoire de Math\'ematiques Jean Leray, UMR 6629 CNRS, Universit\'e de Nantes,\\ 
2, rue de la Houssini\`ere BP 92208, 44322 Nantes - France.\\
E-mail address:petit@math.univ-nantes.fr

\newpage
\section{Introduction}
The main ingredient in the harmonic maps approach of superrigidity for semi-simple Lie groups is the Mok-Siu-Yeung formula \cite{MSY},\cite{SY} for harmonic maps defined on compact locally symmetric spaces of non-compact type (cf. also \cite{JY}). Actually, by means of this formula, it can be shown that a harmonic map from a compact locally symmetric space of non-compact type (with some exceptions) into a Riemannian manifold with nonpositive curvature is rigid in the sense of it is a totally geodesic isometric imbedding. The contact locally sub-symmetric spaces defined by Bieliavsky, Falbel and Gorodski \cite{BFG},\cite{FG},\cite{FGR} are the contact analogues of the riemannian locally symmetric spaces. These spaces can be characterized as contact metric manifolds for which the curvature and the torsion of the Tanaka-Webster connection are parallel in the direction of the contact distribution. Moreover, these spaces are strictly pseudoconvex $CR$ manifolds. In an other hand, in the setting of contact metric manifolds, the analogue of harmonic maps seems to be the pseudoharmonic maps defined by Barletta, Dragomir and Urakawa \cite{BD1},\cite{BDU}. Also the main purpose of this article is to derive Mok-Siu-Yeung type formulas for horizontal maps (i.e. maps preserving the contact distributions) from compact contact locally sub-symmetric spaces into strictly pseudoconvex $CR$ manifolds in order to obtain some rigidity theorems for horizontal pseudoharmonic maps under curvature assumptions. The plan of this article is the following. The section 2 begins to recall basic facts concerning the contact metric manifolds and the strictly pseudoconvex $CR$ manifolds, next, we focus our attention on the pseudo-hermitian curvature tensor of a strictly pseudoconvex $CR$ manifold, which plays a central part in the following. In section 3, we investigate the contact sub-symmetric spaces, the main result of this section (Theorem 3.3) which is related to the work of Cho \cite{JTC}, is an explicit formula for the pseudo-hermitian curvature tensor of a contact locally sub-symmetric space with non zero pseudo-hermitian torsion. The section 4 is devoted to derive Mok-Siu-Yeung type formulas for horizontal maps between strictly pseudoconvex $CR$ manifolds (Proposition 4.2). In section 5, we extend the notion of pseudoharmonic maps defined in \cite{BD1},\cite{BDU} to the setting of horizontal maps between contact metric manifolds and we define the notion of CR-pluriharmonic maps for horizontal maps between strictly pseudoconvex $CR$ manifolds. It is interesting to note that these two notions are strongly related to the Rumin complex \cite{MR}. In section 6, we obtain some rigidity theorems for the horizontal pseudoharmonic maps when the source manifold is a compact contact locally sub-symmetric space. The main result of this section (Theorem 6.1) asserts that any horizontal pseudoharmonic map $\phi$ from a compact contact locally sub-symmetric space of non-compact type, holonomy irreducible and torsionless, (with some exceptions) into a Sasakian manifold with nonpositive pseudo-Hermitian complex sectional curvature satisfies $\nabla d\phi=0$ where $\nabla d\phi$ is the covariant derivative of $d\phi$ with respect to Tanaka-Webster connections. As application (Corollary 6.1) we deduce that $\phi$ preserves some special curves called parabolic geodesics \cite{SD} and therefore $\phi$ is, in some sense, totally geodesic. In section 7, we restrict our attention to CR maps from compact contact locally sub-symmetric spaces into strictly pseudoconvex $CR$ manifolds and we obtain the following rigidity result (Theorem 7.1): any horizontal pseudoharmonic CR map from a compact contact locally sub-symmetric space of non-compact type (with some exceptions) into a pseudo-Hermitian space form with negative pseudo-Hermitian scalar curvature is constant. This article is a first step in the study of horizontal pseudoharmonic maps from compact strictly pseudoconvex $CR$ manifolds into strictly pseudoconvex $CR$ manifolds with nonpositive pseudo-Hermitian sectional curvature. In particular, some existence results are missing for the moment (excepted if the target manifold is Tanaka-Webster flat).\\

The author wants to thanks the CNRS for the d\'el\'egation CNRS that he has benefited during the preparation of this article.
 
 \section{Connection and curvature on contact metric manifolds}

. \underline{Contact metric manifolds}\\

\quad A contact form on a smooth manifold $M$ of dimension $m=2d+1$ is a $1$-form $\theta$ satisfying $\theta\wedge{(d\theta)}^{d}\neq 0$ everywhere on $M$. If $\theta$ is a contact form on $M$, the hyperplan subbundle $H$ of $TM$ given by $H=Ker\,\theta$ is 
called a contact structure. The Reeb field associated to $\theta$ is the unique vector field $\xi$ on $M$ satisfying 
$\theta(\xi)=1$ and $d\theta(\xi,.)=0$. By a contact manifold $(M,\theta)$ we mean a manifold $M$ endowed with a fixed 
contact form $\theta$.\\
If $(M,\theta)$ is a contact manifold then $TM$ decomposes as $TM=H\oplus\R\xi$. Consequently any $p$-tensor $t$ on $M$ decomposes
as $t=t_{H}+t_{\xi}$ with $t_{H}=t\circ\Pi_{H}$ and $t_{\xi}=t\circ\Pi_{\R\xi}$ ($\Pi_{H}$ and $\Pi_{\R\xi}$ are the canonical projections on $H$ and $\R\xi$). The tensors $t_{H}$ and $t_{\xi}$ are respectively called the horizontal part and the vertical part of $t$. 
Note that for an antisymmetric $p$-tensor $\gamma$ we have $\gamma_{\xi}=\theta\wedge i(\xi)\gamma$. We denote by $\wh{*}$ and $\wx{*}$ the bundles of horizontal and vertical antisymmetric tensors and by $\oh{*}$ and $\ox{*}$ the horizontal and vertical forms associated to.\\ 

Let $(M,\theta)$ be a contact manifold, then there exists a riemannian metric $g_{\theta}$ and a $(1,1)$-tensor 
field $J$ on $M$ such that:
$$g_{\theta}(\xi,X)=\theta(X),\quad J^{2}=-Id+\theta\otimes\xi,\quad g_{\theta}(JX,Y)=d\theta(X,Y),\quad X,Y\in TM.$$
The metric $g_{\theta}$ (called the Webster metric) is said to be associated to $\theta$. We call $(\theta,\xi,J,g_{\theta})$ 
a contact metric structure and $(M,\theta,\xi,J,g_{\theta})$ a contact metric manifold (cf. Blair\cite{DB}).
In the following $\omega_{\theta}:=d\theta$.\\
 
Let $(M,\theta,\xi,J,g_{\theta})$ be a contact metric manifold. We define $L:{\Omega}^{k}(M)\to{\Omega}^{k+2}(M)$ by $L=\omega_{\theta}\wedge$. The restriction of $L$ to $\oh{*}$ will be denoted by $L_{H}$ and the adjoint of $L_{H}$ for the usual scalar product on $\oh{*}$ by $\wedge_{H}$.
Recall that, for any $\gamma_{H}\in\oh{p}$,  
$$(\wedge_{H}\gamma_{H})(X_{1},\ldots,X_{p-2})=\frac{1}{2}tr_{H}\gamma_{H}(\,.\, ,J\,.\, ,X_{1},\ldots,X_{p-2}),$$
where $tr_{H}$ is the trace calculated with respect to a $g_{\theta}$-orthonormal frame of $H$.\\

Let ${\Omega}_{H_{0}}^{*}(M):=\{\gamma_{H}\in\oh{*},\ \wedge_{H}\gamma_{H}=0\}$ and ${\Fc}_{\xi}^{*}(M):=\{\gamma_{\xi}\in\ox{*},\ L\gamma_{\xi}=0\}$ be the bundle of primitive horizontal forms on $M$ and the bundle of coprimitive vertical forms on $M$.
We recall the Lefschetz decomposition
$$\oh{*}={\Omega}_{H_{0}}^{*}(M)\oplus L_{H}{\Omega}_{H_{0}}^{*}(M)\oplus\ldots\oplus L_{H}^{d}{\Omega}_{H_{0}}^{*}(M).$$

. \underline{connection and curvature}\\

For the torsion and the curvature of a connection $\nabla$ we adopt the conventions $T(X,Y)=[X,Y]-\nabla_{X}Y+\nabla_{Y}X$ and 
$R(X,Y)=[\nabla_{Y},\nabla_{X}]-\nabla_{[Y,X]}$.\\
In the following, $N$ is the $TM$-valued $2$-form given by: 
$$N(Y,Z)=J^{2}[Y,Z]+[JY,JZ]-J[Y,JZ]-J[JY,Z]+\omega_{\theta}(Y,Z)\xi.$$

\begin{Prop}(Generalized Tanaka-Webster connection cf. \cite{NT},\cite{ST},\cite{SW}). Let $(M,\theta,\xi,J,g_{\theta})$ 
be contact metric manifold, then there exists a unique affine connection $\nabla$ on $TM$ with torsion $T$ 
(called the generalized Tanaka-Webster connection) such that:\\
(a) $\nabla\theta=0$,  $\nabla\xi=0$.\\
(b) $\nabla g_{\theta}=0$.\\
(c) $T_{H}=-\omega_{\theta}\otimes\xi$ and $i(\xi)T=-\frac{1}{2}i(\xi)N$.\\ 
(d) $\displaystyle{(\nabla_{X}\omega_{\theta})(Y,Z)=g_{\theta}((\nabla_{X}J)(Y),Z)=\frac{1}{2}\omega_{\theta}(X,N_{H}(Y,Z))}$ for any $X,Y,Z\in TM$. 
\end{Prop}

The endomorphism $\tau:=i(\xi)T$ is called the generalized Tanaka-Webster torsion or sub-torsion. Note that 
$\tau$ is $g_{\theta}$-symmetric with trace-free and satisfies $\tau\circ J=-J\circ\tau$.\\

A contact metric manifold $(M,\theta,\xi,J,g_{\theta})$ for which $J$ is integrable (i.e. $N_{H}=0$ or equivalently $\nabla J=0$) is 
called a strictly pseudoconvex $CR$ manifold. A strictly pseudoconvex $CR$ manifold for which the Tanaka-Webster torsion vanishes is
called a Sasakian manifold.\\

The curvature $R$ of the generalized Tanaka-Webster connection $\nabla$ satisfies the following Bianchi identities 
(cf. \cite{FG},\cite{NT}):\\ 
(First Bianchi identity) 
\begin{eqnarray}\label{h1}
R_{H}(X,Y)Z+R_{H}(Z,X)Y+R_{H}(Y,Z)X&=&\omega_{\theta}(X,Y)\tau(Z)+\omega_{\theta}(Z,X)\tau(Y)\nonumber\\
                                   &&+\omega_{\theta}(Y,Z)\tau(X)
\end{eqnarray}
\begin{equation}\label{h2}
R(X,\xi)Z+R(\xi,Z)X=(\nabla_{X}\tau)(Z)-(\nabla_{Z}\tau)(X)
\end{equation}
with $(\nabla_{X}\tau)(Z)=\nabla_{X}\tau(Z)-\tau(\nabla_{X}Z)$.\\
 
Remember that any horizontal $2$-tensor $t_{H}$ on $M$ decomposes into $t_{H}=t_{H}^{+}+t_{H}^{-}$, where $t_{H}^{\pm}:=\frac{1}{2}(t_{H}\pm J^{*}t_{H})$ are respectively the $J$-invariant part and the $J$-anti-invariant part of $t_{H}$.  

If $M$ is a strictly pseudoconvex $CR$ manifold, then we have the decomposition:
$$R_{H}=R^{+}_{H}+R^{-}_{H},$$
with 
\begin{equation}\label{h3}
R^{-}_{H}(X,Y)=-\frac{1}{2}\Bigl({(\tau(X))}^{*}\wedge{(JY)}^{*}-{(\tau(Y))}^{*}\wedge{(JX)}^{*}-{(J\tau(X))}^{*}\wedge{Y}^{*}+{(J\tau(Y))}^{*}\wedge{X}^{*}\Bigr).
\end{equation}

Also we define the pseudo-Hermitian curvature tensor $\rw$ of a strictly pseudoconvex $CR$ manifold by:
$$\rw(X,Y,Z,W)=g_{\theta}(R^{+}_{H}(X,Y)Z,W),\quad X,Y,Z,W\in H.$$

In order to give the algebraic properties of $\rw$, we recall some definitions related to the curvature algebra (cf. \cite{JPB}).\\

Let $(V,q)$ be an euclidean space and $Q\in\bigotimes^{4}V^{*}$. We define the Bianchi map $b(Q)$ by:
$$b(Q)(X,Y,Z,W)=Q(X,Y,Z,W)+Q(Z,X,Y,W)+Q(Y,Z,X,W),\quad X,Y,Z,W\in V.$$
Recall that $b(S^{2}({\wedge}^{2}V^{*}))={\wedge}^{4}V^{*}$ and that we have the decomposition:
$$S^{2}({\wedge}^{2}V^{*})=Ker\,b\oplus{\wedge}^{4}V^{*}.$$
Let $h,k\in\bigotimes^{2}V^{*}$. We define the symmetric product $h\odot k\in S^{2}(\bigotimes^{2}V^{*})$ and the Kulkarni product $h\kn k\in\bigotimes^{2}({\wedge}^{2}V^{*})$ respectively by: 
$$(h\odot k)(X,Y,Z,W)=h(X,Y)k(Z,W)+h(Z,W)k(X,Y),$$
and 
$$(h\kn k)(X,Y,Z,W)=(h\odot k)(X,Z,Y,W)-(h\odot k)(X,W,Y,Z).$$
\noindent
Note that if $h,k\in S^{2}V^{*}$ then $h\kn k\in S^{2}({\wedge}^{2}V^{*})$ and $b(h\kn k)=0$.\\ 
If $h,k\in {\wedge}^{2}V^{*}$ then $h\kn k\in S^{2}({\wedge}^{2}V^{*})$ and $b(h\kn k)=-2b(h\odot k)$.\\ 
If $h\in S^{2}V^{*}$ and $k\in {\wedge}^{2}V^{*}$ then $h\kn k\in {\wedge}^{2}({\wedge}^{2}V^{*})$ and $b(h\kn k)=-2b(k\otimes h)$.\\
We define the Ricci contraction $c(Q)$ by:
$$c(Q)(X,Y)=tr\,Q(.,X,.,Y),$$
with the trace is taken with respect to a $q$-orthonormal basis of $V$.\\

Let $(M,\theta,\xi,J,g_{\theta})$ be a contact metric manifold and let $Q_{H}\in\bigotimes^{2}(\wh{2})$. The endomorphism $\qhh$ of $\wh{2}$ associated to $Q_{H}$ is defined by:
$$(\qhh\gamma_{H})(X,Y)=\frac{1}{2}\sum_{i,j}Q_{H}(\ep{i},\ep{j},X,Y)\gamma_{H}(\ep{i},\ep{j}),$$
where $\gamma_{H}\in\wh{2}$, $\{\ep{i}\}$ is a local $g_{\theta}$-orthonormal frame of $H$ and $X,Y\in H$.\\

For $X,Y,Z,W\in H$, we have $\langle\qhh(X\wedge Y),Z\wedge W\rangle=Q_{H}(X,Y,Z,W)$.\\

Let $\who$ be the bundle of primitif horizontal antisymmetric $2$-tensors. For $Q_{H}\in S^{2}(\wh{2})$, we define ${Q}_{H_{0}}\in S^{2}(\who)$ by: 
$${Q}_{H_{0}}=Q_{H}-\frac{1}{d}(\qhh\omega_{\theta})\odot\omega_{\theta}+\frac{1}{2d^2}(\wedge_{H}(\qhh\omega_{\theta}))\omega_{\theta}\odot\omega_{\theta}.$$
Note that, for $Q_{H}$ viewed as horizontal $\wh{2}$-valued $2$-form, we have $\qhh\omega_{\theta}=\wedge_{H}Q_{H}$.\\

Let $\whpm{\pm}$ be the bundle of $J$-invariant ($J$-anti-invariant) horizontal antisymmetric $2$-tensors.
For $Q_{H}\in S^{2}(\wh{2})$, we define $Q_{H}^{\pm}\in S^{2}(\whpm{\pm})$ by:  
\begin{eqnarray*}
Q_{H}^{\pm}(X,Y,Z,W)&=&\frac{1}{4}(Q_{H}(X,Y,Z,W)\pm Q_{H}(JX,JY,Z,W)\pm Q_{H}(X,Y,JZ,JW)\\
&+&Q_{H}(JX,JY,JZ,JW)).
\end{eqnarray*}
Note that 
$${(g_{\theta_{H}}\kn g_{\theta_{H}})}_{0}=g_{\theta_{H}}\kn g_{\theta_{H}}-\frac{1}{d}\omega_{\theta}\odot\omega_{\theta}$$
and 
$${(g_{\theta_{H}}\kn g_{\theta_{H}})}^{\pm}=\frac{1}{2}(g_{\theta_{H}}\kn g_{\theta_{H}}\pm\omega_{\theta}\kn\omega_{\theta}),$$
with  $g_{\theta_{H}}={g_{\theta}}_{/H}$.\\
For any $P_{H},Q_{H}\in S^{2}(\wh{2})$, the scalar product $\langle P_{H},Q_{H}\rangle$ is defined by: 
$$\langle P_{H},Q_{H}\rangle=\frac{1}{2}tr_{H}\widehat{{P}_{H}}\circ\qhh.$$
We have 
$$\langle{(g_{\theta_{H}}\kn g_{\theta_{H}})}^{\pm},Q_{H}\rangle=tr_{H}\widehat{Q_{H}^{\pm}}$$
and for $Q_{H}\in S^{2}(\wh{2})\cap Ker\,b$,
$$\langle\omega_{\theta}\odot\omega_{\theta},Q_{H}\rangle=tr_{H}\widehat{Q_{H}^{+}}-tr_{H}\widehat{Q_{H}^{-}}.$$

If $M$ is $CR$ then $\rw$ can be written by (\ref{h3}) as:
\begin{equation}\label{h4}
\rw=R_{H}+\frac{1}{2}(\omega_{\theta}\kn A_{\theta}-g_{\theta_{H}}\kn B_{\theta}),
\end{equation}
with $R_{H}(X,Y,Z,W)=g_{\theta}(R_{H}(X,Y)Z,W)$, $A_{\theta}(X,Y)=g_{\theta}(\tau(X),Y)$ and $B_{\theta}(X,Y)=\omega_{\theta}(\tau(X),Y)$.\\ 
It follows from (\ref{h1}) and (\ref{h4}) that $\rw$ satisfies the following the algebraic properties:\\

$\rw\in S^{2}(\wh{2})\cap Ker\,b$ and $\rw\in S^{2}(\whpm{+})$.\\

The pseudo-Hermitian Ricci tensor $\riw\in\sh{2}$, the pseudo-Hermitian Ricci form $\row\in\whpm{+}$ and the pseudo-Hermitian scalar curvature $\sw$  are respectively defined by: 
$$\riw=c_{H}(\rw),\quad \row=-\widehat{\rw}\omega_{\theta},\quad \sw=tr_{H}\riw,$$ 
with $c_{H}(Q_{H})$ calculated with respect to a $g_{\theta}$-orthonormal frame of $H$.\\
Now, let $I^{\C}_{H}\in S^{2}(\whpm{+})\cap Ker\,b$ given by: 
$$I^{\C}_{H}=\frac{1}{8}\Bigl(g_{\theta_{H}}\kn g_{\theta_{H}}+\omega_{\theta}\kn\omega_{\theta}+2\omega_{\theta}\odot\omega_{\theta}\Bigr).$$ 

We have the decompositions (cf. \cite{CM}):
$$\rw=\frac{\sw}{d(d+1)}I^{\C}_{H}+\frac{1}{d+2}\Bigl(\frac{1}{2}(\riwo\kn g_{\theta_{H}}-\rowo\kn\omega_{\theta})-\rowo\odot\omega_{\theta}\Bigr)+\Cm$$
$$\rwo=\frac{\sw}{d(d+1)}I^{\C}_{H_{0}}+\frac{1}{2(d+2)}{\Bigl(\riwo\kn g_{\theta_{H}}-\rowo\kn\omega_{\theta}\Bigr)}_{0}+\Cm,$$
where $\riwo$ (respectively $\rowo$) is the traceless part of $\riw$ (respectively the primitive part of $\row$) and 
$\Cm\in S^{2}({\wedge}_{H_{0}}^{2,+}(M))\cap Ker\,b\cap Ker\ c_{H}$.\\

\begin{Rem} The tensor $\Cm$, introduced by Chern and Moser in \cite{CM}, is called the Chern-Moser tensor. Note that $\Cm$
is a pseudo-conformal invariant.
\end{Rem} 

Now, we define the pseudo-Hermitian sectional curvature (resp. pseudo-Hermitian complex sectional curvature) of a $2$-plane $P=\R\{X,Y\}\subset H$ (resp. $P=\C\{Z,W\}\subset H^{\C}$) by: 
$$K^{\Wc}(P)={{\langle\widehat{\rw}(X\wedge Y),X\wedge Y\rangle}\over{\langle X\wedge Y,X\wedge Y\rangle}},\quad K^{\Wc^{\C}}(P)=
{{(\widehat{R^{\Wc^{\C}}_{H}}(Z\wedge W),\overline{Z\wedge W})}\over{(Z\wedge W,\overline{Z\wedge W})}},$$ 
where $(\;,\;)$ and $\widehat{R^{\Wc^{\C}}_{H}}$ are the natural extensions to $\wedge^{2}H^{\C}$ of $\langle\;,\;\rangle$ and $\widehat{\rw}$.\\
The holomorphic pseudo-Hermitian sectional curvature of a holomorphic $2$-plane $P=\R\{X,JX\}\subset H$ is defined by: 
$$HK^{\Wc}(P)={{\langle\widehat{\rw}(X\wedge JX),X\wedge JX\rangle}\over{\langle X\wedge JX,X\wedge JX\rangle}}.$$
We say that a strictly pseudoconvex $CR$ manifold $(M,\theta,\xi,J,g_{\theta})$ has constant holomorphic pseudo-Hermitian sectional curvature if $HK^{\Wc}(P)$ is constant for any holomorphic $2$-plane $P\subset H$ and for any point of $M$. In this case, we have $\displaystyle{\rw=\frac{\sw}{d(d+1)}I^{\C}_{H}}$ with $\sw$ constant (cf. \cite{EB}). Also we call $I^{\C}_{H}$ the holomorphic pseudo-Hermitian curvature tensor.\\

Let $(E,g^{E},\nabla^{E})$ be a riemannian vector bundle over a contact metric manifold $M$ and let ${\Omega}^{*}(M;E)$ (resp $\ohe{*}$) 
be the bundle of $E$-valued forms (resp. horizontal $E$-valued forms) on $M$. We assume that $M$ is endowed with the 
generalized Tanaka-Webster connection $\nabla$. Remember (cf. \cite{RP1}) that for any $\sigma\in{\Omega}^{p}(M;E)$, 
\begin{eqnarray*}
(\nabla_{X}\sigma)(X_{1},\ldots,X_{p})&=&\nabla_{X}^{E}\sigma(X_{1},\ldots,X_{p})-\sum_{i=1}^{p}\sigma(X_{1},\ldots,\nabla_{X}X_{i},\ldots,X_{p})\\
(d^{\nabla^{E}}\sigma)(X_{1},\ldots,X_{p+1})&=&\sum_{i=1}^{p+1}{(-1)}^{i+1}\nabla_{X_{i}}^{E}\sigma(X_{1},\ldots,\hat{X_{i}},\ldots,X_{p+1})\\                                           
                                           &+&\sum_{i<j}{(-1)}^{i+j}\sigma([X_{i},X_{j}],X_{1},\ldots,\hat{X_{i}},\ldots,\hat{X_{j}},\ldots,X_{p+1})\\
                                           &=&\sum_{i=1}^{p+1}{(-1)}^{i+1}(\nabla_{X_{i}}\sigma)(X_{1},\ldots,\hat{X_{i}},\ldots,X_{p+1})\\
                                           &+&\sum_{i<j}{(-1)}^{i+j}\sigma(T(X_{i},X_{j}),X_{1},\ldots,\hat{X_{i}},\ldots,\hat{X_{j}},\ldots,X_{p+1}).
\end{eqnarray*}
Now, for any $\sgh\in\ohe{p}$, we define
\begin{eqnarray*}
(\dhsh)(X_{1},\ldots,X_{p+1})&:=&{(d^{\nabla^{E}}\sgh)}_{H}(X_{1},\ldots,X_{p+1})\\
                             &=&\sum_{i=1}^{p+1}{(-1)}^{i+1}(\nabla_{X_{i}}\sgh)(X_{1},\ldots,\hat{X_{i}},\ldots,X_{p+1})\\
(\lx\sgh)(X_{1},\ldots,X_{p})&:=&(i(\xi)(d^{\nabla^{E}}\sgh))(X_{1},\ldots,X_{p})\\
                             &=&(\nabla_{\xi}\sgh)(X_{1},\ldots,X_{p})+\sum_{i=1}^{p}{(-1)}^{i}\sgh(\tau(X_{i}),X_{1},\ldots,\hat{X_{i}},\ldots,X_{p})\\
(\delta^{\nabla^{E}}_{H}\sgh)(X_{1},\ldots,X_{p-1})&=&-tr_{H}(\nabla_{.\ }\sgh)(\ .\ ,X_{1},\ldots,X_{p-1})\\
(R_{H}(X,Y)\sgh)(X_{1},\ldots,X_{p})&=&R_{H}^{E}(X,Y)\sgh(X_{1},\ldots,X_{p})-\sum_{i=1}^{p}\sgh(X_{1},\ldots,R_{H}(X,Y)X_{i},\ldots,X_{p}),
\end{eqnarray*}
where $R^{E}$ is the curvature of $\nabla^{E}$.\\
Note that for any $\sgh\in\ohe{p}$, we have: 
$${{d}_{H}^{\nabla^{E}}}^{2}\sgh=-L_{H}(\lx\sgh)-R_{H}^{E}\wedge\sgh,$$
where $R_{H}^{E}\wedge\sgh$ is the wedge product of the horizontal $End(E)$-valued $2$-form $R_{H}^{E}$ with the horizontal $E$-valued $p$-form $\sgh$. Now the horizontal $End(H)$-valued $2$-form $R_{H}$ satisfies the second Bianchi identity:
\begin{eqnarray}\label{h5}
(\nabla_{X}R_{H})(Y,Z)+(\nabla_{Z}R_{H})(X,Y)+(\nabla_{Y}R_{H})(Z,X)&=&-\omega_{\theta}(X,Y)(i(\xi)R)(Z)+\omega_{\theta}(Z,X)(i(\xi)R)(Y)\nonumber\\
                                                                     &&-\omega_{\theta}(Y,Z)(i(\xi)R)(X)\nonumber\\
(\nabla_{\xi}R_{H})(X,Y)-(\nabla_{X}i(\xi)R)(Y)+(\nabla_{Y}i(\xi)R)(X)&=&R_{H}(\tau(X),Y)+R_{H}(X,\tau(Y)).
\end{eqnarray}

\begin{Rem} The Bianchi identity (\ref{h5}) is equivalent to $d_{H}^{\nabla}R_{H}=-L_{H}(i(\xi)R)$ and ${\Lc}_{\xi}^{\nabla}R_{H}=d_{H}^{\nabla}(i(\xi)R)$ where $\nabla$ is the connection on $\wh{2}$ induced by the Tanaka-Webster connection on $H$.
\end{Rem}

\section{Contact sub-symmetric spaces}

\begin{Def} A contact (locally) sub-symmetric space is a contact metric manifold $(M,\theta,\xi,J,g_{\theta})$ such that for every point $x_{0}\in M$ there exists an isometry (resp a local isometry) $\psi$, called the sub-symmetry at $x_{0}$, satisfying $\psi(x_{0})=x_{0}$ and $d\psi(x_{0})/H_{x_{0}}=-id$.  
\end{Def}

\begin{Th}\cite{BFG} Let $(M,\theta,\xi,J,g_{\theta})$ be a contact metric manifold endowed with its generalized Tanaka-Webster connection $\nabla$ and let $R$ (respectively $T$) be the curvature (respectively torsion) of $\nabla$. Then:\\
(i) $M$ is a contact locally sub-symmetric space if and only if $\nabla_{H}R=\nabla_{H}T=0$.\\
(ii) If $M$ is a contact locally sub-symmetric space, $\nabla$-complete and simply-connected then $M$ is a contact sub-symmetric space.\\
(iii) If $M$ is a contact sub-symmetric space then $M=G/K$ where $G$ is the closed subgroup of $I(M,g_{\theta})$ generated by all the sub-symmetries $\psi(x_{0})$, $x_{0}\in M$, and $K$ is the isotropy subgroup at a base point (i.e. M is an homogeneous manifold).
\end{Th}

Note that the conditions $\nabla_{H}R=\nabla_{H}T=0$ are equivalent to $\nabla_{H}R_{H}=0$, $i(\xi)R=0$, $\nabla\omega_{\theta}=0$ and $\nabla_{H}\tau=0$. 

\begin{Def}\cite{FGR} A contact sub-symmetric space $M$ is said to be irreducible if the Lie algebra 
$\mathfrak{Hol}(M)$ of the holonomy group $Hol(M)$ acts irreducibly on $H$.
\end{Def}

The simply-connected contact sub-symmetric spaces have been classified by Bieliavsky, Falbel and Gorodski in \cite{BFG}.

\begin{Th} Every simply-connected contact sub-symmetric space of dimension $\geq5$ has the following type:\\
\\
\begin{tabular}{|c|c|c|c|}
\hline
holonomy trivial & torsionless & \multicolumn{2}{c|}{$\Hc_{2p+1}$}\\\hline
& & Compact type  & Non-compact type \\\hline
& & compact Hermitian (CH):  & non-compact Hermitian (NCH):\\
& & $SU(p+q)/SU(p)\times SU(q)$ & $SU(p,q)/SU(p)\times SU(q)$\\
& & $SO(2p)/SU(p)$ & $SO^{*}(2p)/SU(p)$\\ 
holonomy & torsionless & $Sp(p)/SU(p)$ & $Sp(p,\R)/SU(p)$\\ 
irreducible & & $SO(p+2)/SO(p)$\quad$(p\geq 3)$ & $SO_{0}(p,2)/SO(p)$\quad$(p\geq 3)$\\
& & $E_{6(-78)}/Spin(10)$ & $E_{6(-14)}/Spin(10)$\\ 
& & $E_{7(-133)}/E_{6}$ & $E_{7(-25)}/E_{6}$\\\cline{2-4}
& with torsion & $SO(p+1)\ltimes\R^{p+1}/SO(p)$ & $SO_{0}(p,1)\ltimes\R^{p+1}/SO(p)$\\ 
& $(p\geq 3)$ & & $SO_{0}(p+1,1)/SO(p)$\\\hline
holonomy & torsionless & $\Hc_{2p+1}\times_{G}CH$\quad $(G\subset U(p))$ & $\Hc_{2p+1}\times_{G}NCH$\\\cline{2-4}
not & & $SO(4)/SO(2)$ & $SO_{0}(2,2)/SO(2)$\\ 
irreducible & with torsion & $SO(3)\ltimes\R^{3}/SO(2)$ & $SO_{0}(2,1)\ltimes\R^{3}/SO(2)$\\ 
& & & $SO_{0}(3,1)/SO(2)$\\\hline
\end{tabular}
\end{Th}

\begin{Rem} The contact sub-symmetric spaces of compact Hermitian type (respectively non-compact Hermitian type) arise from $S^{1}$-fibrations over irreducible Hermitian symmetric spaces of compact type (respectively non-compact type). The previous list is not complete, because $S^{1}$-fibrations over not irreducible Hermitian symmetric spaces also produce examples of contact sub-symmetric spaces. We do not consider these examples in the following.
\end{Rem}

A contact locally sub-symmetric space $(M,\theta,\xi,J,g_{\theta})$ is always a strictly pseudoconvex $CR$ manifold (since $\nabla\omega_{\theta}=0$). Now we recall the notion of homogeneous strictly pseudoconvex $CR$ manifold and symmetric strictly pseudoconvex $CR$ manifold (cf. \cite{KZ},\cite{EM}).\\

Let $(M,H,J,\theta,g_{\theta})$ and $(N,H^{'},J^{'},\theta^{'},g_{\theta^{'}})$ be strictly pseudoconvex $CR$ manifolds then a map $\phi:M\to N$ such that $d\phi(H)\subset H^{'}$ and $J^{'}\circ\dfh=\dfh\circ J$ is called a CR map from $M$ to $N$ (the definition is also valid in the general context of $CR$ manifolds \cite{AS}). A $CR$ automorphism $\phi:M\to M$ is a diffeomorphism and a CR map from $M$ to $M$. The group of $CR$ automorphisms $Aut_{CR}(M)$ is a Lie group. A $CR$ automorphism $\phi:M\to M$ is called a pseudo-Hermitian transformation if $\phi^{*}\theta=\theta$. The group of pseudo-Hermitian transformations $PsH(M,\theta)$ is a Lie subgroup of $Aut_{CR}(M)$ and also a Lie subgroup of $I(M,g_{\theta})$.

\begin{Def}\cite{EM} A strictly pseudoconvex $CR$ manifold $(M,H,J,\theta,g_{\theta})$ is called homogeneous if there exists a closed subgroup $G$ of $PsH(M,\theta)$ which acts transitively on $M$.
\end{Def}

\begin{Def} A (locally) symmetric strictly pseudoconvex $CR$ manifold is a strictly pseudoconvex $CR$ manifold $(M,H,J,\theta,g_{\theta})$ such that for every point $x_{0}\in M$ there exists a pseudo-Hermitian transformation (resp a local pseudo-Hermitian transformation) $\psi$, called the pseudo-Hermitian symmetry at $x_{0}$, satisfying $\psi(x_{0})=x_{0}$ and $d\psi(x_{0})/H_{x_{0}}=-id$.  
\end{Def}

If $(M,H,J,\theta,g_{\theta})$ is a symmetric strictly pseudoconvex $CR$ manifold then $M=G/K$ where $G$ is the closed subgroup of $PsH(M,\theta)$ generated by all the pseudo-Hermitian symmetries $\psi(x_{0})$, $x_{0}\in M$, and $K$ is the isotropy subgroup at a base point. Also $M$ is an homogeneous strictly pseudoconvex $CR$ manifold. Note that the contact sub-symmetric spaces torsionless are symmetric Sasakian manifolds.\\

Let $(M,\theta,\xi,J,g_{\theta})$ be a simply-connected contact sub-symmetric space and $\Gamma$ be a cocompact discrete subgroup of $PsH(M,\theta)$ acting freely on $M$ then $M/\Gamma$ is a compact contact locally sub-symmetric space.\\

Now we investigate the properties of the pseudo-Hermitian curvature tensor on a contact locally sub-symmetric space.

\begin{Prop} Let $(M,\theta,\xi,J,g_{\theta})$ be a contact locally sub-symmetric space endowed with its 
Tanaka-Webster connection $\nabla$. Then we have $\nabla\rw=0$ and consequently $\sw$ is constant.
\end{Prop}

\noindent{Proof.} Since $\nabla_{H}R_{H}=0$ and $\nabla J=0$, we have $\nabla_{H}R^{+}_{H}=0$. Now, we must prove that $\nabla_{\xi}R^{+}_{H}=0$. Equation (\ref{h5}) together with the assumption $i(\xi)R=0$ gives, for any $X,Y\in H$, 
$$(\nabla_{\xi}R_{H})(X,Y)=R_{H}(\tau(X),Y)+R_{H}(X,\tau(Y)).$$
Since $J\circ\tau=-\tau\circ J$, we deduce that 
\begin{equation}\label{h6}
(\nabla_{\xi}R^{+}_{H})(X,Y)=R^{-}_{H}(\tau(X),Y)+R^{-}_{H}(X,\tau(Y)).
\end{equation}
If $\tau=0$, we have automatically $\nabla_{\xi}R^{+}_{H}=0$. Now, if $\tau\neq 0$, the assumption $\nabla_{H}\tau=0$ implies that $|\tau|^{2}$ is a strictly positive constant and that $\displaystyle{{\tau}^{2}=\frac{|\tau|^{2}}{2d}id_{H}}$ (cf. lemma 1 of \cite{BC}). This assumption together with (\ref{h3}) implies that $R^{-}_{H}(\tau(X),Y)=-R^{-}_{H}(X,\tau(Y))$ and then (\ref{h6}) becomes $\nabla_{\xi}R^{+}_{H}=0$. Hence $\nabla R^{+}_{H}=0$ and $\nabla\rw=0$. $\Box$\\

In \cite{BC}, Boeckx and Cho prove that a contact metric manifold $M$ endowed with its generalized Tanaka-Webster connection $\nabla$ satisfying the conditions $\nabla_{H}J\circ\tau=0$ and $\tau\neq 0$ is a strictly pseudoconvex $CR$ manifold (i.e. $J$ is integrable) and a $(k,\mu)$-space (cf. \cite{BKP}). Moreover, Cho gives, in \cite{JTC}, a formula for the riemannian curvature tensor of $M$ if $M$ has constant holomorphic pseudo-Hermitian sectional curvature. Now we obtain a formula for the pseudo-Hermitian curvature tensor of a strictly pseudoconvex $CR$ manifold satisfying $\nabla_{H}\tau=0$. 

\begin{Th} Let $(M,\theta,\xi,J,g_{\theta})$ be a strictly pseudoconvex $CR$ manifold endowed with its Tanaka-Webster connection $\nabla$. Assume that $\nabla_{H}\tau=0$ and $\tau\neq 0$, then the pseudo-Hermitian curvature tensor and the Chern-Moser tensor of $M$ are given by: 
\begin{equation}\label{h7}
\rw=\frac{\sw}{d^2}\Bigl(I^{\C}_{H}+\frac{2d}{|\tau|^{2}}\Tc_{H}\Bigr),\quad\Cm=\frac{\sw}{d^2}\Bigl(\frac{1}{d+1}I^{\C}_{H_{0}}+\frac{2d}{|\tau|^{2}}\Tc_{H_{0}}\Bigr),
\end{equation}
with $\displaystyle{\Tc_{H}=\frac{1}{8}\Bigl(A_{\theta}\kn A_{\theta}+B_{\theta}\kn B_{\theta}\Bigr)}.$
Moreover, if $d\geq 2$, then $M$ is a contact locally sub-symmetric space.
\end{Th}

The proof of the theorem needs the following Lemma.

\begin{Le} Let $(M,\theta,\xi,J,g_{\theta})$ be a strictly pseudoconvex $CR$ manifold such that $\nabla_{H}\tau=0$
and $\tau\neq 0$. Then $\displaystyle{\row=-\frac{\sw}{2d}\omega_{\theta}}$ (i.e. $M$ is pseudo-Einstein), and 
$\displaystyle{\nabx\tau=-\frac{\sw}{d^2}J\circ\tau}$.
\end{Le}

\noindent{Proof.} First recall that the assumptions $\nabla_{H}\tau=0$ and $\tau\neq 0$ imply that $|\tau|^{2}$ is a strictly positive constant and that $\displaystyle{{\tau}^{2}=\frac{|\tau|^{2}}{2d}id_{H}}$.
Now, we have for any $X,Y\in H$,
$$R_{H}(X,Y)\tau=B_{H}(X,Y)\tau-\omega_{\theta}(X,Y)\nabx\tau,$$ 
with $B_{H}(X,Y)=\nabla_{Y}\nabla_{X}-\nabla_{\nabla_{Y}X}-(\nabla_{X}\nabla_{Y}-\nabla_{\nabla_{X}Y})$.
Since $\nabla_{H}\tau=0$, we have $R_{H}(X,Y)\tau=-\omega_{\theta}(X,Y)\nabx\tau$ and also  
$R^{+}_{H}(X,Y)\tau=-\omega_{\theta}(X,Y)\nabx\tau$. We obtain
$$g_{\theta}(R^{+}_{H}(X,Y)\tau(Z),J\tau(W))-g_{\theta}(\tau(R^{+}_{H}(X,Y)Z),J\tau(W))=-\omega_{\theta}(X,Y)g_{\theta}((\nabx\tau)(Z),J\tau(W)).$$ Hence,
$$\rw(X,Y,\tau(Z),J\tau(W))+\frac{|\tau|^{2}}{2d}\rw(X,Y,Z,JW)=-\omega_{\theta}(X,Y)g_{\theta}((\nabx\tau)(Z),J\tau(W)).$$
Let $\{\ep{i}\}$ be  a local $g_{\theta}$-orthonormal frame of $H$, then 
\begin{equation}\label{h8}
\sum_{1\leq i\leq 2d}\rw(X,Y,\tau(\ep{i}),J\tau(\ep{i}))+\frac{|\tau|^{2}}{2d}\rw(X,Y,\ep{i},J\ep{i})=-\omega_{\theta}(X,Y)\sum_{1\leq i\leq 2d}g_{\theta}((\nabx\tau)(\ep{i}),J\tau(\ep{i})).
\end{equation}
Since $\displaystyle{{\tau}^{2}=\frac{|\tau|^{2}}{2d}id_{H}}$, then (\ref{h8}) becomes 
$$\frac{|\tau|^{2}}{d}\sum_{1\leq i\leq 2d}\rw(\ep{i},J\ep{i},X,Y)=-\omega_{\theta}(X,Y)\sum_{1\leq i\leq 2d}g_{\theta}((\nabx\tau)(\ep{i}),J\tau(\ep{i}))$$ 
which is
$$\row=-\wedge_{H}\rw=\omega_{\theta}\frac{d}{|\tau|^{2}}(\nabx A_{\theta},B_{\theta}).$$
Since $\displaystyle{\wedge_{H}\row=-\frac{\sw}{2}}$, we deduce that $\displaystyle{\row=-\frac{\sw}{2d}\omega_{\theta}}$ and then $M$ is pseudo-Einstein. Now we have 
$$(\wedge_{H}R^{+}_{H})(\tau(X))-\tau((\wedge_{H}R^{+}_{H})(X))=-d(\nabx\tau)(X).$$
The assumption $\displaystyle{\row=-\frac{\sw}{2d}\omega_{\theta}}$ gives
$$(\wedge_{H}R^{+}_{H})(\tau(X))-\tau((\wedge_{H}R^{+}_{H})(X))=\frac{\sw}{d}J\circ\tau(X)$$ 
and then $\displaystyle{\nabx\tau=-\frac{\sw}{d^2}J\circ\tau}$. $\Box$\\

\noindent{Proof of Theorem 3.3} Using (\ref{h2}), we obtain that for any $X,Y,Z\in H$
$$g_{\theta}(R(X,\xi)Y,Z)=g_{\theta}((\nabla_{Z}\tau)(X),Y)-g_{\theta}((\nabla_{Y}\tau)(X),Z).$$ 
By the assumption $\nabla_{H}\tau=0$, we have $i(\xi)R=0$. Now, equation (\ref{h5}) together with $i(\xi)R=0$, $J\circ\tau=-\tau\circ J$ and $\displaystyle{{\tau}^{2}=\frac{|\tau|^{2}}{2d}id_{H}}$ yields 
\begin{equation}\label{h9}
R^{+}_{H}(\tau(X),\tau(Y))+\frac{|\tau|^{2}}{2d}R^{+}_{H}(X,Y)=(\nabx R^{-}_{H})(X,\tau(Y)).
\end{equation}
By (\ref{h3}), we have 
\begin{eqnarray*}
(\nabx R^{-}_{H})(X,Y)&=&-\frac{1}{2}\Bigl({((\nabx\tau)(X))}^{*}\wedge{(JY)}^{*}-{((\nabx\tau)(Y))}^{*}\wedge{(JX)}^{*}\\
&-&{((J(\nabx\tau)(X))}^{*}\wedge{Y}^{*}+{((J(\nabx\tau)(Y))}^{*}\wedge{X}^{*}\Bigr).
\end{eqnarray*}
We have $\displaystyle{\nabx\tau=-\frac{\sw}{d^2}J\circ\tau}$ (Lemma 3.1), it follows that
$$(\nabx R^{-}_{H})(X,\tau(Y))=\frac{\sw}{4d^3}|\tau|^{2}\Bigl({X}^{*}\wedge{Y}^{*}+{(JX)}^{*}\wedge{(JY)}^{*}\Bigr)
+\frac{\sw}{2d^2}\Bigl({(\tau(X))}^{*}\wedge{(\tau(Y))}^{*}+{(J\tau(X))}^{*}\wedge{(J\tau(Y))}^{*}\Bigr).$$
Then (\ref{h9}) becomes
\begin{eqnarray}\label{h10}
\rw(\tau(X),\tau(Y),Z,W)+\frac{|\tau|^{2}}{2d}\rw(X,Y,Z,W)&=&\frac{\sw}{8d^3}|\tau|^{2}\Bigl(g_{\theta_{H}}\kn g_{\theta_{H}}+\omega_{\theta}\kn\omega_{\theta}\Bigr)(X,Y,Z,W)\nonumber\\
                                                &+&\frac{\sw}{4d^2}\Bigl(A_{\theta}\kn A_{\theta}+B_{\theta}\kn B_{\theta}\Bigr)(X,Y,Z,W).
\end{eqnarray}
Now, we have
\begin{eqnarray*}
g_{\theta}(R^{+}_{H}(X,Y)\tau(Z),\tau(W))-g_{\theta}(\tau(R^{+}_{H}(X,Y)Z),\tau(W))&=&-\omega_{\theta}(X,Y)g_{\theta}((\nabx\tau)(Z),\tau(W))\\
                                                                  &=&-\frac{\sw}{2d^3}|\tau|^{2}\omega_{\theta}(X,Y)\omega_{\theta}(Z,W).
\end{eqnarray*}
Hence
\begin{equation}\label{h11}
\rw(\tau(X),\tau(Y),Z,W)-\frac{|\tau|^{2}}{2d}\rw(X,Y,Z,W)=-\frac{\sw}{4d^3}|\tau|^{2}(\omega_{\theta}\odot\omega_{\theta})(X,Y,Z,W).
\end{equation}
We deduce from (\ref{h10}) and (\ref{h11}) the following expression for the pseudo-Hermitian curvature
\begin{eqnarray*}
\rw(X,Y,Z,W)&=&\frac{\sw}{8d^2}\Bigl(g_{\theta_{H}}\kn g_{\theta_{H}}+\omega_{\theta}\kn\omega_{\theta}+2\omega_{\theta}\odot\omega_{\theta}\Bigr)(X,Y,Z,W)\\
&+&\frac{\sw}{4d|\tau|^{2}}\Bigl(A_{\theta}\kn A_{\theta}+B_{\theta}\kn B_{\theta}\Bigr)(X,Y,Z,W).
\end{eqnarray*}
The expression for $\Cm$ directly follows from the decomposition of $\rw$. Now we assume $d>1$. Since $\nabla_{H}\tau=0$, we have the formula $\delta_{H}\riw=-\frac{1}{2}d_{H}\sw$. Also $M$ pseudo-Einstein and $d>1$ yields to $\sw$ constant. Since $\nabla_{H}g_{\theta}=\nabla_{H}\omega_{\theta}=\nabla_{H}A_{\theta}=\nabla_{H}B_{\theta}=0$ and $\sw$ is constant, then we have by the previous formula for $\rw$ and (\ref{h4}) that $\nabla_{H}R_{H}=0$. Consequently $M$ is a contact locally sub-symmetric space. $\Box$

\begin{Cor} The pseudo-Hermitian curvature tensor $\rw$ of a contact locally sub-symmetric space $M$ has the following form. If $M$ is holonomy irreducible and torsionless, in this case $M$ is the total space of a $S^{1}$-fibration $\pi$ over an irreducible Hermitian locally symmetric space $B$ and $\rw$ is given by $\rw=\pi^{*}R^{B}$ where $R^{B}$ is the curvature of $B$. If $M$ has torsion, in this case $\rw$ is given by formula (\ref{h7}). Note that, in each case, $M$ is pseudo-Einstein with $\sw$ constant. 
\end{Cor}

\begin{Rem} Note that the contact sub-symmetric spaces of non-compact Hermitian type have nonpositive pseudo-Hermitian sectional curvature whereas the contact sub-symmetric spaces of compact Hermitian type have nonnegative pseudo-Hermitian sectional curvature.
\end{Rem}

\section{Mok-Siu-Yeung type formulas for horizontal maps between strictly pseudoconvex $CR$ manifolds}

Let $(M,\theta,\xi,J,g_{\theta},\nabla)$ be a contact metric manifold endowed with the (generalized) Tanaka-Webster 
connection and let $(E,g^{E},\nabla^{E})$ be a Riemannian vector bundle over $M$.
For any $Q_{H}\in\wh{2}\otimes End(E)$ and $\sgh\in\wh{1}\otimes E$, we define $Q_{H}(\sgh)\in\wh{1}\otimes E$ by:
$$Q_{H}(\sgh)(X)=\sum_{i}Q_{H}(\ep{i},X)\sgh(\ep{i}),$$
where $\{\ep{i}\}$ is a local $g_{\theta}$-orthonormal frame of $H$.
For any $Q_{H}\in S^{2}(\wh{2})$ and $s_{H}\in\sh{2}\otimes E$ (respectively $\sgh\in\wh{2}\otimes E$), we define 
$\qhc s_{H}\in\sh{2}\otimes E$ (respectively $\qhh\sgh\in\wh{2}\otimes E$) by:
\begin{eqnarray*}
(\qhc s_{H})(X,Y)&=&\sum_{i,j}Q_{H}(\ep{i},X,Y,\ep{j})s_{H}(\ep{i},\ep{j}),\\
(\qhh\sgh)(X,Y)&=&\frac{1}{2}\sum_{i,j}Q_{H}(\ep{i},\ep{j},X,Y)\sgh(\ep{i},\ep{j}).\\
\end{eqnarray*}

\begin{Prop} Let $(M,\theta,\xi,J,g_{\theta},\nabla)$ be a compact contact metric manifold and let $Q_{H}\in\Gamma({S}^{2}(\wh{2}))$ satisfying as horizontal $\wh{2}$-valued $2$-form the assumptions $\delta^{\nabla}_{H}Q_H=0$ and $\wedge_{H}Q_{H}=0$. Then, for any $\sgh\in{\Omega}_{H}^{1}(M;E)$, we have:
\begin{eqnarray}\label{h12}
\int_{M}\langle\qhc\nbsh,\nbsh\rangle+\langle\widehat{(b(Q_{H})-Q_{H})}\dhsh,\dhsh\rangle v_{g_{\theta}}&=&2\int_{M}\langle (\qhh R_{H}^{E})(\sgh),\sgh\rangle\nonumber\\
&&-\langle (c_{H}(\widehat{R_{H}}\circ\qhh){)}^{\Sc},\sgh^{*}g^{E}\rangle v_{g_{\theta}},\nonumber\\
\end{eqnarray}
where, for any horizontal $2$-tensor $\mu_{H}$, ${\mu}_{H}^{\Sc}(X,Y)=\mu_{H}(X,Y)+\mu_{H}(Y,X)$. 
\end{Prop}

\noindent{Proof.} Let $Q_{H}\in\Gamma({S}^{2}(\wh{2}))$ satisfying $\delta^{\nabla}_{H}Q_H=0$ and $\wedge_{H}Q_{H}=0$, then formula (8) of \cite{RP2} gives for any $\sgh\in{\Omega}_{H}^{1}(M;E)$
$$(\delta_{H}Q_{H}(\nabla))\sgh=\Rc^{Q}_{H}\sgh,$$
where $Q_{H}(\nabla)$ and $\Rc^{Q}_{H}$ are given in a local orthonormal frame $\{\ep{i}\}$ of $H$ by:
$$Q_{H}(\nabla)(X)=\sum_{i}\qhh(\ep{i}\wedge X).\nab{i}=\frac{1}{2}\sum_{i,k,l}Q_{H}(\ep{i},X,\ep{k},\ep{l})\ep{k}.\ep{l}.\nab{i},$$
and
$$\Rc^{Q}_{H}=-\frac{1}{2}\sum_{i,j}\qhh(\ep{i}\wedge\ep{j}).R_{H}(\ep{i},\ep{j})=-\frac{1}{4}\sum_{i,j,k,l}\qh{i}{j}{k}{l}\ep{k}.\ep{l}.R_{H}(\ep{i},\ep{j}).$$
By integrating, we obtain
$$\int_{M}\langle Q_{H}(\nabla)\sgh,\nabla\sgh\rangle v_{g_{\theta}}=\int_{M}\langle\Rc^{Q}_{H}\sgh,\sgh\rangle v_{g_{\theta}}.$$
We have
$$\langle Q_{H}(\nabla)\sgh,\nabla\sgh\rangle=-\frac{1}{2}\sum_{i,j,k,l}\qh{i}{j}{k}{l}\langle\ep{l}.\nab{i}\sgh,\ep{k}.\nab{j}\sgh\rangle.$$
Now, for any $X,Y\in TM$ and any $\sigma,\gamma\in{\Omega}^{1}(M;E)$, we have (cf. \cite{RP1}):
\begin{equation}\label{h13}
\langle X.\sigma,Y.\gamma\rangle=g_{\theta}(X,Y)\langle\sigma,\gamma\rangle+\langle i(X)\sigma,i(Y)\gamma\rangle-\langle i(Y)\sigma,i(X)\gamma\rangle.
\end{equation}
We deduce from (\ref{h13}) that 
\begin{equation}\label{h14}
\langle Q_{H}(\nabla)\sgh,\nabla\sgh\rangle=-\sum_{i,j,k,l}\qh{i}{j}{k}{l}\langle(\nab{i}\sgh)(\ep{l}),(\nab{j}\sgh)(\ep{k})\rangle
\end{equation}
Now we have for any $\sgh\in{\Omega}_{H}^{1}(M;E)$ and any $X,Y\in H$
$$(\nabla_{X}\sgh)(Y)=\frac{1}{2}\Bigl((\nbsh)(X,Y)+(\dhsh)(X,Y)\Bigr),$$
with $(\nbsh)(X,Y)=(\nabla_{X}\sgh)(Y)+(\nabla_{Y}\sgh)(X)$ and $(\dhsh)(X,Y)=(\nabla_{X}\sgh)(Y)-(\nabla_{X}\sgh)(Y)$. 
Then (\ref{h14}) becomes
\begin{eqnarray*}
\langle Q_{H}(\nabla)\sgh,\nabla\sgh\rangle&=&-\frac{1}{4}\Bigl(\sum_{i,j,k,l}\qh{i}{j}{k}{l}\langle\nbish{i}{l},\nbish{j}{k}\rangle\\
&&+\frac{1}{2}\sum_{i,j,k,l}(\qh{i}{j}{k}{l}+\qh{k}{i}{j}{l})\langle\dhish{i}{l},\dhish{j}{k}\rangle\\
&&+2\sum_{i,j,k,l}\qh{i}{j}{k}{l}\langle\nbish{i}{l},\dhish{j}{k}\rangle\Bigr)\\
                                            &=&-\frac{1}{4}\Bigl(\sum_{i,j,k,l}\qh{i}{j}{k}{l}\langle\nbish{i}{l},\nbish{j}{k}\rangle\\
&&+\frac{1}{2}\sum_{i,j,k,l}(b({Q}_{H})(\ep{i},\ep{l},\ep{j},\ep{k})-\qh{i}{l}{j}{k})\langle\dhish{i}{l},\dhish{j}{k}\rangle\\
&&+2\sum_{i,j,k,l}\qh{i}{j}{k}{l}\langle\nbish{i}{l},\dhish{j}{k}\rangle\Bigr)\\
                                            &=&-\frac{1}{4}\Bigl(\sum_{j,k}\langle(\qhc\nbsh)(\ep{j},\ep{k}),\nbish{j}{k}\rangle\\
&&+\sum_{j,k}\langle\widehat{(b(Q_{H})-Q_{H})}\dhsh)(\ep{j},\ep{k}),\dhish{j}{k}\rangle\\
&&+2\sum_{j,k}\langle(\qhc\nbsh)(\ep{j},\ep{k}),\dhish{j}{k}\rangle\Bigr)\\
                                            &=&-\frac{1}{2}\Bigl(\langle\qhc\nbsh,\nbsh\rangle+\langle\widehat{(b(Q_{H})-Q_{H})}\dhsh,\dhsh\rangle\Bigr).
\end{eqnarray*}
For the second term, we have 
\begin{eqnarray*}
\langle\Rc^{Q}_{H}\sgh,\sgh\rangle&=&\frac{1}{4}\sum_{i,j,k,l}\qh{i}{j}{k}{l}\langle\ep{l}.R_{H}(\ep{i},\ep{j})\sgh,\ep{k}.\sgh\rangle\\
                                  &=&-\frac{1}{2}\sum_{i,j,k,l}\qh{i}{j}{k}{l}\langle (R_{H}(\ep{i},\ep{j})\sgh)(\ep{k}),\sgh(\ep{l})\rangle.
\end{eqnarray*}
Since $(R_{H}(X,Y)\sgh)(Z)=R_{H}^{E}(X,Y)\sgh(Z)-\sgh(R_{H}(X,Y)Z)$, we have 
\begin{eqnarray*}
\langle\Rc^{Q}_{H}\sgh,\sgh\rangle&=&-\frac{1}{2}\sum_{i,j,k,l}\qh{i}{j}{k}{l}\langle R_{H}^{E}(\ep{i},\ep{j})\sgh(\ep{k}),\sgh(\ep{l})\rangle\\
&&+\frac{1}{2}\sum_{i,j,k,l,m}\qh{i}{j}{k}{l}R_{H}(\ep{i},\ep{j},\ep{k},\ep{m})\langle\sgh(\ep{l}),\sgh(\ep{m})\rangle\\
&=&-\sum_{k,l}\langle(\qhh R_{H}^{E})(\ep{k},\ep{l})\sgh(\ep{k}),\sgh(\ep{l})\rangle\\
&&+\frac{1}{2}\sum_{l,m}\Bigl(c_{H}(\widehat{R_{H}}\circ\qhh)(\ep{l},\ep{m})+c_{H}(\widehat{R_{H}}\circ\qhh)(\ep{m},\ep{l})\Bigr)(\sgh^{*}g^{E})(\ep{l},\ep{m})\\
&=&-\langle (\qhh R_{H}^{E})(\sgh),\sgh\rangle+\langle (c_{H}(\widehat{R_{H}}\circ\qhh){)}^{\Sc},\sgh^{*}g^{E}\rangle.
\end{eqnarray*}
Hence the formula. $\Box$\\

Assume that $(M,\theta,\xi,J,g_{\theta},\nabla)$ and $(N,\theta^{'},\xi^{'},J^{'},g_{\theta^{'}},\nabla^{'})$ 
are contact metric manifolds endowed with their Tanaka-Webster connections and let $\phi:M\to N$ be a differential map.
Let $\ftn$ be the pull-back bundle of $TN$ endowed with the metric and the connection induced by those of $TN$. The 
covariant derivative of the $\ftn$-valued $1$-form $\dfh$ is given by:
$$(\nabla_{X}\dfh)(Y)={\nabla}^{{'}^{\ftn}}_{X}\dfh(Y)-\dfh(\nabla_{X}Y),$$
where ${\nabla}^{{'}^{\ftn}}$ denotes the connection induced by ${\nabla}^{'}$ on $\ftn$.

\begin{Def} Let $H=Ker\,\theta$ and $H^{'}=Ker\,\theta^{'}$. A map $\phi:M\to N$ such that $d\phi(H)\subset H^{'}$  
is called a horizontal map from $M$ to $N$. We denote by $\Hc(M,N)$ the subspace of horizontal maps from $M$ to $N$.  
\end{Def}
Note that a horizontal map $\phi$ satisfies $\phi^{*}\theta^{'}=f\theta$ with $f\in C^{\infty}(M,\R)$.

\begin{Le} For any horizontal map $\phi:M\to N$, we have:
\begin{equation}\label{h15}
\dhfh=-\omega_{\theta}\otimes\dfxh,
\end{equation}
and
\begin{equation}\label{h16}
{\Lc}_{\xi}^{\nabla^{'}}\dfh={\nabla}_{H}^{{'}^{\ftn}}\dfxh-f\,\tau^{'}\circ\dfh.
\end{equation}
\end{Le}

\noindent{Proof.} For any map $\phi:M\to N$, we have ${d}^{\nabla^{'}}d\phi=-\phi^{*}T^{'}$ where $T^{'}$ is the torsion of $\nabla^{'}$. Hence
\begin{eqnarray*}
{d}^{\nabla^{'}}\dfh&=&{d}^{\nabla^{'}}(d\phi-\theta\otimes d\phi(\xi))=-\phi^{*}T^{'}-\omega_{\theta}\otimes d\phi(\xi)+\theta\wedge{\nabla}^{{'}^{\ftn}}d\phi(\xi)\\
                    &=&-{(\phi^{*}T^{'})}_{H}-\omega_{\theta}\otimes d\phi(\xi)+\theta\wedge({\nabla}^{{'}^{\ftn}}d\phi(\xi)-i(\xi)\phi^{*}T^{'}).
\end{eqnarray*}
We deduce that 
$$\dhfh={({d}^{\nabla^{'}}\dfh)}_{H}=-{(\phi^{*}T^{'})}_{H}-\omega_{\theta}\otimes d\phi(\xi)$$
and 
$${\Lc}_{\xi}^{\nabla^{'}}\dfh=i(\xi)({d}^{\nabla^{'}}\dfh)={\nabla}_{H}^{{'}^{\ftn}}d\phi(\xi)-i(\xi)\phi^{*}T^{'}.$$
Now, let $\phi:M\to N$ be a horizontal map, then we have $\phi^{*}\theta^{'}=f\theta$ and $\phi^{*}\omega_{\theta^{'}}=f\omega_{\theta}-\theta\wedge df_{H}$. Since $T^{'}=-\omega_{\theta^{'}}\otimes\xi^{'}+\theta^{'}\wedge\tau^{'}$, we deduce that 
$\phi^{*}T^{'}=-f\omega_{\theta}\otimes\xi^{'}+\theta\wedge(df_{H}\otimes\xi^{'}+f\,\tau^{'}\circ\dfh)$ 
and consequently
$$\dhfh=\omega_{\theta}\otimes(f\xi^{'}-d\phi(\xi))=-\omega_{\theta}\otimes\dfxh.$$
Now,
\begin{eqnarray*}
{\Lc}_{\xi}^{\nabla^{'}}\dfh&=&{\nabla}_{H}^{{'}^{\ftn}}d\phi(\xi)-df_{H}\otimes\xi^{'}-f\,\tau^{'}\circ\dfh\\
                            &=&{\nabla}_{H}^{{'}^{\ftn}}(d\phi(\xi)-f\xi^{'})-f\,\tau^{'}\circ\dfh={\nabla}_{H}^{{'}^{\ftn}}\dfxh-f\,\tau^{'}\circ\dfh.\ \Box
\end{eqnarray*}
\\
In the following, for any horizontal symmetric $2$-tensor $\mu_{H}$, we denote by ${\mu}_{H_{0}}$ its traceless part.
\newpage
\begin{Prop}(Mok-Siu-Yeung type formulas for horizontal maps between strictly pseudoconvex $CR$ manifolds) 
Let $(M,\theta,\xi,J,g_{\theta},\nabla)$ and $(N,\theta^{'},\xi^{'},J^{'},g_{\theta^{'}},\nabla^{'})$ be strictly pseudoconvex $CR$ manifolds with the assumption $M$ compact. For any $Q_{H_{0}}^{+}\in\Gamma({S}^{2}({\wedge}_{H_{0}}^{2,+}(M)))$ (resp. $Q_{H}^{-}\in\Gamma({S}^{2}(\whpm{-}))$) satisfying $\delta^{\nabla}_{H}Q_{H_{0}}^{+}=0$ and ${(c_{H}(Q_{H_{0}}^{+}))}_{0}=0$ (resp. $\delta^{\nabla}_{H}Q_{H}^{-}=0$ and ${(c_{H}(Q_{H}^{-}))}_{0}=0$) and any horizontal map $\phi$ from $M$ to $N$, we have:
\begin{eqnarray}\label{h17}
&&\int_{M}\langle\stackrel{\circ}{Q_{H_{0}}^{+}}{(\nbfh)}_{0},{(\nbfh)}_{0}\rangle-\frac{tr_{H}\widehat{Q_{H_{0}}^{+}}}{d^2}\Bigl(\nn{\delfh}+d^2\nn{\dfxh}\Bigr) v_{g_{\theta}}\nonumber\\
&&=4\int_{M}2\langle Q_{H_{0}}^{+},\rwpf\rangle-\frac{1}{2}\langle(c_{H}(\widehat{\rw}\circ\widehat{Q_{H_{0}}^{+}}){)}^{\Sc},{(\phi^{*}g_{\theta^{'}})}_{H}\rangle-\langle\stackrel{\circ}{Q_{H_{0}}^{+}}{(\phi^{*}B_{\theta^{'}})}_{H},{(\phi^{*}g_{\theta^{'}})}_{H}\rangle v_{g_{\theta}},\nonumber\\ \\
&&\int_{M}\langle\stackrel{\circ}{Q_{H}^{-}}{(\nbfh)}_{0},{(\nbfh)}_{0}\rangle-\frac{tr_{H}\widehat{Q_{H}^{-}}}{d^2}\Bigl(\nn{\delfh}-d^2\nn{\dfxh}\Bigr) v_{g_{\theta}}\nonumber\\
&&=4\int_{M}2\langle Q_{H}^{-},\rwpf\rangle-\langle\stackrel{\circ}{Q_{H}^{-}}({(\phi^{*}B_{\theta^{'}})}_{H}-B_{\theta})+\frac{tr_{H}\widehat{Q_{H}^{-}}}{d}B_{\theta},{(\phi^{*}g_{\theta^{'}})}_{H}\rangle v_{g_{\theta}}.\label{h18}
\end{eqnarray}
\end{Prop}

\noindent{Proof.} Let $\phi:M\to N$ be a horizontal map. For any $Q_{H}\in\Gamma({S}^{2}(\wh{2}))$, we obtain, using the relation $\displaystyle{\nbfh={(\nbfh)}_{0}-\frac{1}{d}g_{\theta_{H}}\otimes\delfh}$, that 
\begin{eqnarray*}
\langle\stackrel{\circ}{Q_{H}}\nbfh,\nbfh\rangle&=&\langle\stackrel{\circ}{Q_{H}}{(\nbfh)}_{0},{(\nbfh)}_{0}\rangle+\frac{2}{d}\langle c_{H}(Q_{H})\otimes\delfh,{(\nbfh)}_{0}\rangle\\
                                                &-&\frac{tr_{H}\widehat{Q_{H}}}{d^2}\nn{\delfh}.
\end{eqnarray*}
The assumption ${(c_{H}(Q_{H}))}_{0}=0$ implies that $\langle c_{H}(Q_{H})\otimes\delfh,{(\nbfh)}_{0}\rangle=0$.
Hence we have 
\begin{equation}\label{h19}
\langle\stackrel{\circ}{Q_{H}}\nbfh,\nbfh\rangle=\langle\stackrel{\circ}{Q_{H}}{(\nbfh)}_{0},{(\nbfh)}_{0}\rangle-\frac{tr_{H}\widehat{Q_{H}}}{d^2}\nn{\delfh}.
\end{equation}
Since $\phi:M\to N$ is horizontal then (\ref{h15}) yields $\dhfh=-\omega_{\theta}\otimes\dfxh$. 
Now we have
\begin{eqnarray*}
\langle\widehat{(b(Q_{H})-Q_{H})}\dhfh,\dhfh\rangle&=&\langle\widehat{(b(Q_{H})-Q_{H})}\omega_{\theta},\omega_{\theta}\rangle\nn{\dfxh}\\
&=&\langle b(Q_{H})-Q_{H},\omega_{\theta}\odot\omega_{\theta}\rangle\nn{\dfxh}\\
&=&(tr_{H}\widehat{Q_{H}^{-}}-tr_{H}\widehat{Q_{H}^{+}})\nn{\dfxh}.
\end{eqnarray*}
For $Q_{H_{0}}^{+}\in\Gamma({S}^{2}({\wedge}_{H_{0}}^{2,+}(M)))$ and $Q_{H}^{-}\in\Gamma({S}^{2}(\whpm{-}))$, we have 
\begin{equation}\label{h20}
\langle\widehat{(b(Q_{H_{0}}^{+})-Q_{H_{0}}^{+})}\dhfh,\dhfh\rangle=-tr_{H}\widehat{Q_{H_{0}}^{+}} 
\end{equation}
and 
\begin{equation}\label{h21}
\langle\widehat{(b(Q_{H}^{-}-Q_{H}^{-})}\dhfh,\dhfh\rangle=tr_{H}\widehat{Q_{H}^{-}}.
\end{equation}
Now for $Q_{H}\in\Gamma({S}^{2}(\wh{2}))$, we have 
$\langle (\qhh R_{H}^{\ftn})(\dfh),\dfh\rangle=2\langle Q_{H},{(\phi^{*}\rhp)}_{H}^{\Sc}\rangle$. 
Since $\displaystyle{\rhp=\rwp-\frac{1}{2}(\omega_{\theta^{'}}\kn A_{\theta^{'}}-g_{\theta^{'}}\kn B_{\theta^{'}})}$, we have 
$${(\phi^{*}\rhp)}_{H}=\rwpf+\frac{1}{2}{(\phi^{*}g_{\theta^{'}})}_{H}\kn {(\phi^{*}B_{\theta^{'}})}_{H}-\frac{1}{2}{(\phi^{*}\omega_{\theta^{'}})}_{H}\kn {(\phi^{*}A_{\theta^{'}})}_{H}.$$ 
Since ${(\phi^{*}\omega_{\theta^{'}})}_{H}\kn {(\phi^{*}A_{\theta^{'}})}_{H}\in\Gamma({\wedge}^{2}(\wh{2}))$, we have
\begin{eqnarray}\label{h22}
\langle (\qhh R_{H}^{\ftn})(\dfh),\dfh\rangle&=&4\langle Q_{H},\rwpf\rangle+2\langle Q_{H},{(\phi^{*}g_{\theta^{'}})}_{H}\kn {(\phi^{*}B_{\theta^{'}})}_{H}\rangle\nonumber\\
&=&4\langle Q_{H},\rwpf\rangle-2\langle\stackrel{\circ}{Q_{H}}{(\phi^{*}B_{\theta^{'}})}_{H},{(\phi^{*}g_{\theta^{'}})}_{H}\rangle\end{eqnarray}
The relations $\widehat{g_{\theta_{H}}\kn B_{\theta}}\circ\widehat{\omega_{\theta}\kn\omega_{\theta}}=2\widehat{\omega_{\theta}\kn A_{\theta}}$ and $\widehat{g_{\theta_{H}}\kn B_{\theta}}\circ\widehat{g_{\theta_{H}}\kn g_{\theta_{H}}}=2\widehat{g_{\theta_{H}}\kn B_{\theta}}$ yield to  
$$\widehat{R_{H}}=\widehat{\rw}+\frac{1}{2}\widehat{g_{\theta_{H}}\kn B_{\theta}}\circ\widehat{{(g_{\theta_{H}}\kn g_{\theta_{H}})}^{-}}.$$ Since for any $T_{H}^{\pm},Q_{H}^{\pm}\in S^{2}(\whpm{\pm})$, $\widehat{T_{H}^{\pm}}\circ\widehat{Q_{H}^{\mp}}=0$ and $\widehat{{(g_{\theta_{H}}\kn g_{\theta_{H}})}^{-}}\circ\widehat{Q_{H}^{-}}=2\widehat{Q_{H}^{-}}$ then
\begin{equation}\label{h23}
\langle (c_{H}(\widehat{R_{H}}\circ\widehat{Q_{H_{0}}^{+}}){)}^{\Sc},{(\phi^{*}g_{\theta^{'}})}_{H}\rangle=\langle (c_{H}(\widehat{\rw}\circ\widehat{Q_{H_{0}}^{+}}){)}^{\Sc},{(\phi^{*}g_{\theta^{'}})}_{H}\rangle,
\end{equation}
and 
$$\langle (c_{H}(\widehat{R_{H}}\circ\widehat{Q_{H}^{-}}){)}^{\Sc},{(\phi^{*}g_{\theta^{'}})}_{H}\rangle=\langle (c_{H}(\widehat{g_{\theta_{H}}\kn B_{\theta}}\circ\widehat{Q_{H}^{-}}){)}^{\Sc},{(\phi^{*}g_{\theta^{'}})}_{H}\rangle.$$
Now, we have
$$\langle (c_{H}(\widehat{g_{\theta_{H}}\kn B_{\theta}}\circ\widehat{Q_{H}^{-}}){)}^{\Sc},{(\phi^{*}g_{\theta^{'}})}_{H}\rangle
=\langle(\stackrel{\circ}{c_{H}(Q_{H}^{-})}\circ\stackrel{\circ}{B_{\theta}}{)}^{\Sc}-2\stackrel{\circ}{Q_{H}^{-}}B_{\theta},{(\phi^{*}g_{\theta^{'}})}_{H}\rangle,$$
where $\stackrel{\circ}{\mu_{H}}$ is the symmetric endomorphism associated by $g_{\theta_{H}}$ to the symmetric $2$-tensor $\mu_{H}$. Using the assumption ${(c_{H}(Q_{H}^{-}))}_{0}=0$, we deduce that
\begin{equation}\label{h24}
\langle (c_{H}(\widehat{R_{H}}\circ\widehat{Q_{H}^{-}}){)}^{\Sc},{(\phi^{*}g_{\theta^{'}})}_{H}\rangle=2\langle\frac{tr_{H}\widehat{Q_{H}^{-}}}{d}B_{\theta}-\stackrel{\circ}{Q_{H}^{-}}B_{\theta},{(\phi^{*}g_{\theta^{'}})}_{H}\rangle.
\end{equation}
We obtain the formulas by replacing (\ref{h19}),(\ref{h20}),(\ref{h21}),(\ref{h22}),(\ref{h23}) and (\ref{h24}) in (\ref{h12}). $\Box$\\

As applications of formulas (\ref{h17}) and (\ref{h18}), we recover, in a different way, the Siu formula given in \cite{RP1} and we derive a Tanaka-Weitzenbock formula. First we have the following lemma whose proof is left to the reader.

\begin{Le} We have for any $s_{H}\in\sh{2}\otimes E$ the relations:
\begin{eqnarray*}
\stackrel{\circ}{g_{\theta_{H}}\kn g_{\theta_{H}}}s_{H}&=&2(s_{H}-g_{\theta_{H}}\otimes tr_{H}s_{H}),\quad\stackrel{\circ}{\omega_{\theta}\kn\omega_{\theta}}s_{H}=\stackrel{\circ}{\omega_{\theta}\odot\omega_{\theta}}s_{H}=-2(s_{H}^{+}-s_{H}^{-}),\\
c_{H}(g_{\theta_{H}}\kn g_{\theta_{H}})&=&2(2d-1)g_{\theta_{H}},\quad c_{H}(\omega_{\theta}\kn\omega_{\theta})=c_{H}(\omega_{\theta}\odot\omega_{\theta})=2g_{\theta_{H}},\\
tr_{H}\,\widehat{g_{\theta_{H}}\kn g_{\theta_{H}}}&=&2d(2d-1),\quad tr_{H}\,\widehat{\omega_{\theta}\kn\omega_{\theta}}=tr_{H}\,\widehat{\omega_{\theta}\odot\omega_{\theta}}=2d.
\end{eqnarray*}
\end{Le}

 \begin{Prop} For any horizontal map $\phi$ from a compact Sasakian manifold $M$ to a Sasakian manifold $N$, we have: 
\begin{eqnarray}\label{h25}
\int_{M}\nn{{(\nbfh)}_{0}^{+}}-\Bigl(1-\frac{1}{d}\Bigr)\nn{\delfh}+d(d-1)\nn{\dfxh}v_{g_{\theta}}=4\int_{M}\rf{2}{0}v_{g_{\theta}},\\
\int_{M}\frac{1}{d}\nn{{(\nbfh)}_{0}^{+}}+\Bigl(1-\frac{1}{d}\Bigr)\nn{{(\nbfh)}^{-}}-\Bigl(1-\frac{1}{d^2}\Bigr)\nn{\delfh}-(d^2-1)\nn{\dfxh}v_{g_{\theta}}\nonumber\\
=4\int_{M}\frac{1}{d}\rf{2}{0}+\Bigl(1-\frac{1}{d}\Bigr)\rf{1}{1}-\frac{1}{2}\Bigl(1-\frac{1}{d}\Bigr)\langle\dfh\circ\stackrel{\circ}{\riw},\dfh\rangle v_{g_{\theta}},\label{h26}
\end{eqnarray}
where, in an adapted frame $\{\ep{1},\ldots\ep{d},J\ep{1},\ldots J\ep{d}\}$ of $H$,
$$\rf{2}{0}=\sum_{i,j\leq d}(\rwpfhc(\Zi\wedge\Zj),\overline{\Zi\wedge\Zj})$$
and 
$$\rf{1}{1}=\sum_{i,j\leq d}(\rwpfhc(\Zi\wedge\overline{\Zj}),\overline{\Zi\wedge\overline{\Zj}}),$$
with $\Zi={{1}\over{\sqrt{2}}}(\epsilon_{i}-\sqrt{-1}J\epsilon_{i})$.
\end{Prop}

\noindent{Proof.} Let $Q_{H}^{-}\in\Gamma({S}^{2}(\whpm{-}))$ and $Q_{H_{0}}^{+}\in\Gamma({S}^{2}({\wedge}_{H_{0}}^{2,+}(M)))$ defined by: $$Q_{H}^{-}={(g_{\theta_{H}}\kn g_{\theta_{H}})}^{-}=\frac{1}{2}(g_{\theta_{H}}\kn g_{\theta_{H}}-\omega_{\theta}\kn\omega_{\theta})$$ 
and 
$$Q_{H_{0}}^{+}={(g_{\theta_{H}}\kn g_{\theta_{H}})}_{0}^{+}=\frac{1}{2}(g_{\theta_{H}}\kn g_{\theta_{H}}+\omega_{\theta}\kn\omega_{\theta}-\frac{2}{d}\omega_{\theta}\odot\omega_{\theta}).$$ 
Then we have ${\nabla}_{H}Q_{H_{0}}^{+}={\nabla}_{H}Q_{H}^{-}=0$ (since $\nabla g_{\theta_{H}}=\nabla\omega_{\theta}=0$). Moreover, using Lemma 4.2, we have $c_{H}(Q_{H}^{-})=2(d-1)g_{\theta_{H}}$ and $\displaystyle{c_{H}(Q_{H_{0}}^{+})=2d\Bigl(1-\frac{1}{d^2}\Bigr)g_{\theta_{H}}}$. Also ${(c_{H}(Q_{H}^{-}))}_{0}={(c_{H}(Q_{H_{0}}^{+}))}_{0}=0$. Moreover, we have $\displaystyle{(c_{H}(\widehat{\rw}\circ\widehat{Q_{H_{0}}^{+}}){)}^{\Sc}=4\Bigl(1-\frac{1}{d}\Bigr)\riw}$. Let $\phi$ be a horizontal map from $M$ to $N$, by Lemma 4.2, we have 
\begin{eqnarray*}
\langle\stackrel{\circ}{Q_{H}^{-}}{(\nbfh)}_{0},{(\nbfh)}_{0}\rangle&=&2\nn{{(\nbfh)}_{0}^{+}}\\
\langle\stackrel{\circ}{Q_{H_{0}}^{+}}{(\nbfh)}_{0},{(\nbfh)}_{0}\rangle&=&2\Bigl(\frac{1}{d}\nn{{(\nbfh)}_{0}^{+}}+\Bigl(1-\frac{1}{d}\Bigr)\nn{{(\nbfh)}^{-}}\Bigr)\\
tr_{H}\widehat{Q_{H}^{-}}&=&2d(d-1),\quad tr_{H}\widehat{Q_{H_{0}}^{+}}=2\Bigl(d^2-1\Bigr)\\
\langle(c_{H}(\widehat{\rw}\circ\widehat{Q_{H_{0}}^{+}}){)}^{\Sc},{(\phi^{*}g_{\theta^{'}})}_{H}\rangle&=&2\Bigl(1-\frac{1}{d}\Bigr)\langle\dfh\circ\stackrel{\circ}{\riw},\dfh\rangle\\
\langle Q_{H}^{-},\rwpf\rangle&=&\langle{(g_{\theta_{H}}\kn g_{\theta_{H}})}^{-},\rwpf\rangle=tr_{H}\rwpfh{-}\\
\langle Q_{H_{0}}^{+},\rwpf\rangle&=&\langle{(g_{\theta_{H}}\kn g_{\theta_{H}})}^{+}-\frac{1}{d}\omega_{\theta}\odot\omega_{\theta},\rwpf\rangle\\
&=&\Bigl(1-\frac{1}{d}\Bigr)tr_{H}\rwpfh{+}+\frac{1}{d}tr_{H}\rwpfh{-}.
\end{eqnarray*}
In an adapted frame $\{\ep{1},\ldots,\ep{d},J\ep{1},\ldots,J\ep{d}\}$ of $H$, we have 
$$tr_{H}\rwpfh{\pm}=\sum_{i,j\leq d}({(\phi^{*}\rwp)}_{H}^{\pm}(\ep{i},\ep{j},\ep{i},\ep{j})+{(\phi^{*}\rwp)}_{H}^{\pm}(\ep{i},J\ep{j},\ep{i},J\ep{j})).$$
Now we have, for any $T_{H}\in S^{2}(\wh{2})\cap Ker\,b$, the relations 
$$(\widehat{T_{H}}^{\C}(Z\wedge W),\overline{Z\wedge W})=T_{H}^{-}(X,Y,X,Y)+T_{H}^{-}(X,JY,X,JY)$$
and 
$$(\widehat{T_{H}}^{\C}(Z\wedge\overline{W}),\overline{Z\wedge\overline{W}})=T_{H}^{+}(X,Y,X,Y)+T_{H}^{+}(X,JY,X,JY),$$
with $Z={{1}\over{\sqrt{2}}}(X-\sqrt{-1}JX),\quad W={{1}\over{\sqrt{2}}}(Y-\sqrt{-1}JY)$.
Since $\rwpf\in S^{2}(\wh{2})\cap Ker\,b$, we deduce that
$$\langle Q_{H}^{-},\rwpf\rangle=\sum_{i,j\leq d}(\rwpfhc(\Zi\wedge\Zj),\overline{\Zi\wedge\Zj})$$
and
\begin{eqnarray*}
\langle Q_{H_{0}}^{+},\rwpf\rangle&=&\sum_{i,j\leq d}\Bigl(\frac{1}{d}(\rwpfhc(\Zi\wedge\Zj),\overline{\Zi\wedge\Zj})\\
&+&\Bigl(1-\frac{1}{d}\Bigr)(\rwpfhc(\Zi\wedge\overline{\Zj}),\overline{\Zi\wedge\overline{\Zj}})\Bigr), 
\end{eqnarray*}
with $\Zi={{1}\over{\sqrt{2}}}(\epsilon_{i}-\sqrt{-1}J\epsilon_{i})$.
By replacing in (\ref{h17}) and (\ref{h18}) together with the assumptions $M,N$ Sasakian yields the formulas. $\Box$
 
\section{Horizontal pseudoharmonic maps, CR-pluriharmonic maps and Rumin complex}

. \underline{Pseudoharmonic maps}\\

In \cite{BD1} and \cite{BDU} Barletta, Dragomir and Urakawa have introduced the notion of pseudoharmonic maps from a compact contact metric manifold into a Riemannian manifold. Now we extend this notion to horizontal maps between contact metric manifolds.\\
  
Assume that $(M,\theta,\xi,J,g_{\theta},\nabla)$ and $(N,\theta^{'},\xi^{'},J^{'},g_{\theta^{'}},\nabla^{'})$ are contact metric manifolds endowed with their Tanaka-Webster connections and that $M$ is compact. For any differential map $\phi:M\to N$, 
we define $\dfhh(X)={(\dfh(X))}_{H^{'}}$ with $X\in H$ and the horizontal energy $E_{H,H^{'}}(\phi)$ by: 
$$E_{H,H^{'}}(\phi)=\frac{1}{2}\int_{M}\nn{\dfhh}v_{g_{\theta}}.$$ 

\begin{Prop} For any variation $\phi_{t}$ of $\phi$, we have:
$$\frac{d}{dt}{E_{H,H^{'}}(\phi_{t})}_{|t=0}=\int_{M}g_{\theta^{'}}(\delta^{\nabla^{'}}_{H}\dfhh+i(\phi^{*}\theta^{'})\tau^{'}\circ\dfhh-\xi^{'}tr_{H}{(\phi^{*}A_{\theta^{'}})}_{H},v)v_{g_{\theta}},$$
with $\displaystyle{v={\frac{\partial\phi}{\partial t}}_{|t=0}}$.
\end{Prop}

\noindent{Proof.} Let $\{\phi_{t}\}_{|t|<\epsilon}$ be a variation of $\phi$. We consider the map $\Phi:]-\epsilon,\epsilon[\times M\to N$
given by $\Phi(t,x)=\phi_{t}(x)$ and the pull-back bundle $\Ftn\to]-\epsilon,\epsilon[\times M$ of $TN$ by $\Phi$. Let $\{\ep{i}\}$ 
be a local $g_{\theta}$-orthonormal frame of $H$, then we have
$$\frac{1}{2}\frac{d}{dt}\nn{{d\phi_{t}}_{H,H^{'}}}=\frac{1}{2}\frac{\partial}{\partial t}\sum_{i}g_{\theta^{'}}({(d\Phi(\ep{i}))}_{H^{'}},{(d\Phi(\ep{i}))}_{H^{'}})=\sum_{i}g_{\theta^{'}}({\nabla}^{{'}^{\Ftn}}_{\frac{\partial}{\partial t}}{(d\Phi(\ep{i}))}_{H^{'}},{(d\Phi(\ep{i}))}_{H^{'}}).$$
We have
$$({d}^{\nabla^{'}}d\Phi)(\frac{\partial}{\partial t},\ep{i})={\nabla}^{{'}^{\Ftn}}_{\frac{\partial}{\partial t}}d\Phi(\ep{i})-{\nabla}^{{'}^{\Ftn}}_{\ep{i}}d\Phi(\frac{\partial}{\partial t})-d\Phi([\frac{\partial}{\partial t},\ep{i}])=-T^{'}(d\Phi(\frac{\partial}{\partial t}),d\Phi(\ep{i})),$$
where $T^{'}$ is the torsion of $\nabla^{'}$. Since $[\frac{\partial}{\partial t},\ep{i}]=0$ and ${\nabla}^{{'}^{\Ftn}}$ preserves $H^{'}$, we obtain
$${\nabla}^{{'}^{\Ftn}}_{\frac{\partial}{\partial t}}{(d\Phi(\ep{i}))}_{H^{'}}-{\nabla}^{{'}^{\Ftn}}_{\ep{i}}{(d\Phi(\frac{\partial}{\partial t}))}_{H^{'}}=-{(T^{'}(d\Phi(\frac{\partial}{\partial t}),d\Phi(\ep{i})))}_{H^{'}}.$$
Using $T^{'}=-\omega_{\theta^{'}}\otimes\xi^{'}+\theta^{'}\wedge\tau^{'}$, we obtain
\begin{eqnarray*}
\frac{1}{2}\frac{d}{dt}\nn{{d\phi_{t}}_{H,H^{'}}}&=&\sum_{i}g_{\theta^{'}}\Bigl({\nabla}^{{'}^{\Ftn}}_{\ep{i}}{(d\Phi(\frac{\partial}{\partial t}))}_{H^{'}}-\theta^{'}(d\Phi(\frac{\partial}{\partial t}))\tau^{'}(d\Phi(\ep{i}))\\
&+&\theta^{'}(d\Phi(\ep{i}))\tau^{'}(d\Phi(\frac{\partial}{\partial t})),{(d\Phi(\ep{i}))}_{H^{'}}\Bigr)\\
                                                 &=&\sum_{i}\Bigl(\ep{i}g_{\theta^{'}}({(d\Phi(\frac{\partial}{\partial t}))}_{H^{'}},{(d\Phi(\ep{i}))}_{H^{'}})-g_{\theta^{'}}({(d\Phi(\frac{\partial}{\partial t}))}_{H^{'}},{(d\Phi(\nab{i}\ep{i}))}_{H^{'}})\Bigr)\\
                                                 &-&\sum_{i}g_{\theta^{'}}({(d\Phi(\frac{\partial}{\partial t}))}_{H^{'}},{\nabla}^{{'}^{\Ftn}}_{\ep{i}}{(d\Phi(\ep{i}))}_{H^{'}}-{(d\Phi(\nab{i}\ep{i}))}_{H^{'}})\\
                                                 &-&(\Phi^{*}\theta^{'})(\frac{\partial}{\partial t})\sum_{i}A_{\theta^{'}}({(d\Phi(\ep{i}))}_{H^{'}},{(d\Phi(\ep{i}))}_{H^{'}})\\
                                                 &+&\sum_{i}(\Phi^{*}\theta^{'})(\ep{i})A_{\theta^{'}}({(d\Phi(\frac{\partial}{\partial t}))}_{H^{'}},{(d\Phi(\ep{i}))}_{H^{'}})\\
                                                 &=&-\delta_{H}\alpha_{{(d\Phi(\frac{\partial}{\partial t}))}_{H^{'}}}-g_{\theta^{'}}\Bigl({(d\Phi(\frac{\partial}{\partial t}))}_{H^{'}},\sum_{i}{(({\nabla}^{{'}^{\Ftn}}_{\ep{i}}d\Phi)(\ep{i}))}_{H^{'}}\Bigr)\\
                                                 &-&(\Phi^{*}\theta^{'})(\frac{\partial}{\partial t})\sum_{i}A_{\theta^{'}}({(d\Phi(\ep{i}))}_{H^{'}},{(d\Phi(\ep{i}))}_{H^{'}})\\
                                                 &+&\sum_{i}(\Phi^{*}\theta^{'})(\ep{i})A_{\theta^{'}}({(d\Phi(\frac{\partial}{\partial t}))}_{H^{'}},{(d\Phi(\ep{i}))}_{H^{'}}),
\end{eqnarray*}
with $\alpha_{{(d\Phi(\frac{\partial}{\partial t}))}_{H^{'}}}(X)=g_{\theta^{'}}({(d\Phi(\frac{\partial}{\partial t}))}_{H^{'}},{(d\Phi(X))}_{H^{'}})$. 
We deduce that
\begin{eqnarray*}
\frac{1}{2}\frac{d}{dt}{\nn{{d\phi_{t}}_{H,H^{'}}}}_{|t=0}&=&-\delta_{H}\alpha_{v_{H^{'}}}-g_{\theta^{'}}\Bigl(\sum_{i}{(({\nabla}^{{'}^{\ftn}}_{\ep{i}}d\phi)(\ep{i}))}_{H^{'}},v_{H^{'}}\Bigr)\\
                                                          &&-\theta^{'}(v)\sum_{i}A_{\theta^{'}}(\dfhh(\ep{i}),\dfhh(\ep{i}))+\sum_{i}(\phi^{*}\theta^{'})(\ep{i})A_{\theta^{'}}(\dfhh(\ep{i}),v_{H^{'}})\\
                                                          &=&-\delta_{H}\alpha_{v_{H^{'}}}+g_{\theta^{'}}(\delta^{\nabla^{'}}_{H}\dfhh+i(\phi^{*}\theta^{'})\tau^{'}\circ\dfhh-\xi^{'}tr_{H}{(\phi^{*}A_{\theta^{'}})}_{H},v).
\end{eqnarray*}
The result follows by integrating. $\Box$

\begin{Def} A map $\phi:M\to N$ is called a pseudoharmonic map if it is a critical point of $E_{H,H^{'}}$.
\end{Def}
A map $\phi:M\to N$ is pseudoharmonic if and only if 
$$\delta^{\nabla^{'}}_{H}\dfhh+i(\phi^{*}\theta^{'})\tau^{'}\circ\dfhh=0\quad{\rm and}\quad tr_{H}{(\phi^{*}A_{\theta^{'}})}_{H}=0.$$
A horizontal map $\phi:M\to N$ is pseudoharmonic if and only if $\delta^{\nabla^{'}}_{H}\dfh=0$ and\\$tr_{H}{(\phi^{*}A_{\theta^{'}})}_{H}=0$.\\

. \underline{CR-pluriharmonic maps}\\

Let $(M,\theta,\xi,J,g_{\theta})$ be a strictly pseudoconvex $CR$ manifold of dimension $2d+1$. A real function $h$ on $M$ is called a CR-pluriharmonic function if $h$ is the real part of a CR function on $M$.\\ 

We have the following equivalent characterizations for the CR-pluriharmonic functions.  

\begin{Th}(Lee\cite{JL})\\
The following assertions are equivalent:\\
(i) $h$ is CR-pluriharmonic.\\ 
(ii) There exists a real function $\lambda$ such that $d(J^{*}{dh}_{H}+\lambda\theta)=0$.\\ 
(iii) ${(d_{H}J^{*}{dh}_{H})}_{0}=0$ and $\displaystyle{({\Lc}_{\xi}+\frac{1}{d}d_{H}\delta_{H,J})J^{*}{dh}_{H}={\Lc}_{\xi}J^{*}{dh}_{H}+\frac{1}{d}d_{H}\delta_{H}{dh}_{H}=0}$,
\end{Th}
where ${(d_{H}J^{*}{dh}_{H})}_{0}$ is the primitive part of $d_{H}J^{*}{dh}_{H}$ and $\delta_{H,J}=[\wedge_{H},d_{H}]$.\\

Note that, if $d>1$, then the assumption ${(d_{H}J^{*}{dh}_{H})}_{0}=0$ implies that $\displaystyle{{\Lc}_{\xi}J^{*}{dh}_{H}+\frac{1}{d}d_{H}\delta_{H}{dh}_{H}=0}$ and, if $d=1$, then the assumption ${(d_{H}J^{*}{dh}_{H})}_{0}=0$ is always satisfied for any $h$. 

\begin{Def} Let $(M,\theta,\xi,J,g_{\theta},\nabla)$ and $(N,\theta^{'},\xi^{'},J^{'},g_{\theta^{'}},\nabla^{'})$ be strictly pseudoconvex $CR$ manifolds endowed with their Tanaka-Webster connections together with $dim\,M>3$. A horizontal map $\phi:M\to N$ such that ${({d}_{H}^{\nabla^{'}}J^{*}\dfh)}_{0}=0$ is called a CR-pluriharmonic map from $M$ to $N$.
\end{Def}

\begin{Prop}(i) A map $\phi:M\to N$ is CR-pluriharmonic if and only if 
$${(\nbfh)}_{0}^{+}:={(\nbfh)}^{+}+\frac{1}{d}g_{\theta_{H}}\otimes\delta^{\nabla^{'}}_{H}\dfh=0.$$
(ii) Any CR-pluriharmonic map $\phi:M\to N$ satisfies:
$${\Lc}_{\xi}^{\nabla^{'}}J^{*}\dfh+\frac{1}{d}{d}_{H}^{\nabla^{'}}\delfh=\frac{2}{d-1}tr_{H}\,{(R_{H}^{\ftn})}^{-}(\quad,\,.\,)\dfh(\,.\,)-f\,\tau^{'}\circ J^{*}\dfh.$$
(iii) Any CR map $\phi:M\to N$ is CR-pluriharmonic and we have  
\begin{equation}\label{h27}
\delta_{H}^{\nabla^{'}}J^{*}\dfh=J^{'}\delfh=d\dfxh.
\end{equation}
\end{Prop}

\noindent{Proof.} Recall that for $\gamma_{H}\in\oh{2}$, its primitive part $\gamma_{H_{0}}\in{\Omega}_{H_{0}}^{2}(M)$ is given by
$\displaystyle{\gamma_{H_{0}}=\gamma_{H}-\frac{1}{d}L_{H}\wedge_{H}\gamma_{H}}$.
We have
$${({d}_{H}^{\nabla^{'}}J^{*}\dfh)}_{0}={d}_{H}^{\nabla^{'}}J^{*}\dfh-\frac{1}{d}L_{H}\wedge_{H}{d}_{H}^{\nabla^{'}}J^{*}\dfh
={d}_{H}^{\nabla^{'}}J^{*}\dfh-\frac{1}{d}L_{H}\delhj J^{*}\dfh.$$
Now since $\delhj J^{*}\dfh=\delhj(\dfh\circ J)=\delfh$, we deduce that 
\begin{equation}\label{h28}
{({d}_{H}^{\nabla^{'}}J^{*}\dfh)}_{0}={d}_{H}^{\nabla^{'}}J^{*}\dfh-\frac{1}{d}\omega_{\theta}\otimes\delfh. 
\end{equation}
Now, using (\ref{h15}), we have 
\begin{eqnarray*}
({d}_{H}^{\nabla^{'}}J^{*}\dfh)(JX,Y)&=&(\nabla_{JX}\dfh)(JY)+(\nabla_{Y}\dfh)(X)=\frac{1}{2}\Bigl((\nbfh)(X,Y)+(\nbfh)(JX,JY)\Bigr)\\
                                     &=&{(\nbfh)}^{+}(X,Y).
\end{eqnarray*}
Hence by (\ref{h28}) we obtain
$${({d}_{H}^{\nabla^{'}}J^{*}\dfh)}_{0}(JX,Y)={(\nbfh)}^{+}(X,Y)+\frac{1}{d}g_{\theta_{H}}(X,Y)\delta^{\nabla^{'}}_{H}\dfh.$$
The assumption (i) directly follows. Now we have 
$$\delhj({({d}_{H}^{\nabla^{'}}J^{*}\dfh)}_{0})=\wedge_{H}{d}_{H}^{\nabla^{'}}({({d}_{H}^{\nabla^{'}}J^{*}\dfh)}_{0})
=\wedge_{H}{{d}_{H}^{\nabla^{'}}}^{2}J^{*}\dfh-\frac{1}{d}\wedge_{H}{d}_{H}^{\nabla^{'}}L_{H}\delfh.$$
Since ${{d}_{H}^{\nabla^{'}}}^{2}J^{*}\dfh=-L_{H}(\lx J^{*}\dfh)-R_{H}^{\ftn}\wedge J^{*}\dfh$ and $[{d}_{H}^{\nabla^{'}},L_{H}]=0$, 
we obtain 
\begin{eqnarray}\label{h29}
\delhj({({d}_{H}^{\nabla^{'}}J^{*}\dfh)}_{0})&=&-\wedge_{H}(R_{H}^{\ftn}\wedge J^{*}\dfh)-\wedge_{H}L_{H}(\lx J^{*}\dfh)-\frac{1}{d}\wedge_{H}L_{H}{d}_{H}^{\nabla^{'}}\delfh\nonumber\\
&=&-\wedge_{H}(R_{H}^{\ftn}\wedge J^{*}\dfh)-(d-1)(\lx J^{*}\dfh+\frac{1}{d}{d}_{H}^{\nabla^{'}}\delfh).
\end{eqnarray}
Let $\{\ep{i}\}$ be a local $g_{\theta}$-orthonormal frame of $H$, then we have 
$$\wedge_{H}(R_{H}^{\ftn}\wedge J^{*}\dfh)=\frac{1}{2}\sum_{i}\Bigl(\rhp(\dfh(\ep{i}),\dfh(J\ep{i}))\dfh(JX)+2\rhp(\dfh(\ep{i}),\dfh(X))\dfh(\ep{i})\Bigr)$$
Using (\ref{h1}), we obtain that
$$\sum_{i}\rhp(\dfh(\ep{i},)\dfh(J\ep{i}))\dfh(JX)=2\sum_{i}\rhp(\dfh(JX),\dfh(J\ep{i}))\dfh(\ep{i})+2f(d-1)\tau^{'}((\dfh\circ J)(X)).$$
Hence we have
\begin{eqnarray}\label{h30}
\wedge_{H}(R_{H}^{\ftn}\wedge J^{*}\dfh)&=&-\sum_{i}\Bigl(\rhp(\dfh(X),\dfh(\ep{i}))\dfh(\ep{i})-\rhp(\dfh(JX),\dfh(J\ep{i}))\dfh(\ep{i})\Bigr)\nonumber\\
&+&f(d-1)\tau^{'}((\dfh\circ J)(X))\nonumber\\
&=&-2\sum_{i}{(R_{H}^{\ftn})}^{-}(X,\ep{i})\dfh(\ep{i})+f(d-1)(\tau^{'}\circ J^{*}\dfh)(X).
\end{eqnarray}
The assumption $\phi$ CR-pluriharmonic together with (\ref{h29}) and (\ref{h30}) gives the formula. Hence (ii). 
Let $\phi:M\to N$ be a CR map then $J^{'}\circ\dfh=\dfh\circ J$. Consequently, we have 
$$J^{'}\circ{({d}_{H}^{\nabla^{'}}\dfh)}_{0}={({d}_{H}^{\nabla^{'}}(\dfh\circ J))}_{0}$$
and  
$$J^{'}\delta_{H}^{\nabla^{'}}\dfh=\delta_{H}^{\nabla^{'}}(\dfh\circ J).$$ 
By (\ref{h15}) we have ${({d}_{H}^{\nabla^{'}}\dfh)}_{0}=0$ and we obtain that ${({d}_{H}^{\nabla^{'}}(\dfh\circ J)}_{0}=0$. Now we have 
$$\delta_{H}^{\nabla^{'}}(\dfh\circ J)=-\delhj(\dfh)=-\wedge_{H}({d}_{H}^{\nabla^{'}}\dfh)=d\dfxh.\quad\Box$$
 
The following theorem holds for pseudoharmonic maps between Sasakian manifolds.

\begin{Th} Let $M$ and $N$ be Sasakian manifolds. Assume that $M$ is compact and $N$ has nonpositive pseudo-Hermitian complex sectional curvature. Then:\\
(i) Any horizontal pseudoharmonic map $\phi$ from $M$ to $N$ is CR-pluriharmonic.\\
(ii) If the pseudo-Hermitian Ricci tensor of $M$ is nonnegative, then any horizontal pseudoharmonic map $\phi$ from $M$ to $N$ satisfies $\nabla d\phi=0$ and $|d\phi|=const$.\\
(iii) If the pseudo-Hermitian Ricci tensor of $M$ is positive, then any horizontal pseudoharmonic map $\phi$ from $M$ to $N$ satisfies $\dfh=0$ and consequently $rg_{x}(\phi)\leq 1$, where $rg_{x}(\phi)$ is the rank of $\phi$ at a point $x$ of $M$.
\end{Th}

\noindent{Proof.} Let $\phi$ be a horizontal pseudoharmonic map from $M$ to $N$. If $N$ has nonpositive pseudo-Hermitian complex sectional curvature then (\ref{h25}) yields to $\dfxh=0$ and ${(\nbfh)}_{0}^{+}=0$. In particular, $\phi$ is CR-pluriharmonic. Moreover, if $\riw$ is nonnegative, it follows from (\ref{h26}) that ${(\nbfh)}^{-}=0$. Consequently $\nbfh=0$ and $\dhfh=0$. We deduce that $\nabla_{H}\dfh=0$. The assumptions $M$ and $N$ torsionless together with $\dfxh=0$, yield by (\ref{h16}) that $\nabx\dfh=0$. Taking into account that $d\phi(\xi)=f\xi^{'}$, we obtain that $i(\xi)\phi^{*}\omega_{\theta^{'}}=-df_{H}=0$ and so $f$ is constant. We immediately deduce that $\nabla d\phi=\nabla\dfh+\theta\otimes\nabla d\phi(\xi)=0$ and so $|d\phi|=const$. If $\riw$ is positive, it directly follows from $\langle\dfh\circ\stackrel{\circ}{\riw},\dfh\rangle=0$ that $\dfh=0$. Since $d\phi(\xi)=f\xi^{'}$, we deduce that the rank of $\phi$ is less than equal $1$ at each point of $M$. $\Box$\\

. \underline{Horizontal maps and twisted Rumin pseudo-complex}\\

Let $M$ be a strictly pseudoconvex $CR$ manifold of dimension $2d+1$. We recall that the Rumin complex \cite{MR} is the complex:
$$0\rightarrow\R\rightarrow C^{\infty}(M)\stackrel{{d}_{\Rc}}{\rightarrow}\Rc^{1}(M)\stackrel{{d}_{\Rc}}{\rightarrow}\ldots\stackrel{{d}_{\Rc}}{\rightarrow}\Rc^{d}(M)\stackrel{{D}_{\Rc}}{\rightarrow}\Rc^{d+1}(M)\stackrel{{d}_{\Rc}}{\rightarrow}\ldots\stackrel{{d}_{\Rc}}{\rightarrow}\Rc^{2d+1}(M)\rightarrow 0,$$
where 
\begin{eqnarray*}
\Rc^{p}(M)&=&{\Omega}_{H_{0}}^{p}(M)\quad{\rm for}\ p\leq d\\
          &=&{\Fc}_{\xi}^{p}(M)\quad{\rm for}\ p\geq d+1,
\end{eqnarray*}
and
\begin{eqnarray*}
\begin{array}{ccc}
{d}_{\Rc}\gamma_{H}={({d}_{H}\gamma_{H})}_{0}=({d}_{H}-\frac{1}{d-p+1}L_{H}\delta_{H,J})\gamma_{H}  & {\rm for}\ \gamma_{H}\in\Rc^{p}(M) & (p\leq d-1)\\
{D}_{\Rc}\gamma_{H}=\theta\wedge({\Lc}_{\xi}+d_{H}\delta_{H,J})\gamma_{H} & {\rm for}\ \gamma_{H}\in\Rc^{d}(M) &\\
{d}_{\Rc}\gamma_{\xi}=\theta\wedge i(\xi)(d\gamma_{\xi}) & {\rm for}\ \gamma_{\xi}\in\Rc^{p}(M) & (p\geq d+1).
\end{array}
\end{eqnarray*}   
The formal adjoints of ${d}_{\Rc}$ and ${D}_{\Rc}$ for the usual scalar product are denoted by ${\delta}_{\Rc}$ and ${D}_{\Rc}^{*}$.
The laplacians associated to this complex are defined by:  
\begin{eqnarray*}
\begin{array}{ccc}
{\triangle}_{\Rc}=(d-p){d}_{\Rc}{\delta}_{\Rc}+(d-p+1){\delta}_{\Rc}{d}_{\Rc} & {\rm on}\ \Rc^{p}(M) & (p\leq d-1)\\
{\triangle}_{\Rc}={D}_{\Rc}^{*}{D}_{\Rc}+{({d}_{\Rc}{\delta}_{\Rc})}^{2} & {\rm on}\ \Rc^{d}(M) &\\
{\triangle}_{\Rc}={D}_{\Rc}{D}_{\Rc}^{*}+{({\delta}_{\Rc}{d}_{\Rc})}^{2} & {\rm on}\ \Rc^{d+1}(M) &\\
{\triangle}_{\Rc}=(d-p+1){d}_{\Rc}{\delta}_{\Rc}+(d-p){\delta}_{\Rc}{d}_{\Rc} & {\rm on}\ \Rc^{p}(M) & (p\geq d+2).
\end{array}
\end{eqnarray*}   
The fondamental fact is that, if $M$ is compact, then (cf. \cite{MR}):
$$H_{dR}^{*}(M,\R)=H_{\Rc}^{*}(M,\R)=Ker\,{\triangle}_{\Rc},$$ 
where $H_{dR}^{*}(M,\R)$ and $H_{\Rc}^{*}(M,\R)$ are respectively the cohomologies of the De Rham complex and the Rumin complex.\\
 
Let $(E,\nabla^{E})$ be a vector bundle over $M$ then the previous definitions of $\Rc^{p}(M)$, ${d}_{\Rc}$ and ${D}_{\Rc}$ can be extended to $E$-twisted bundles. Also we define the sequence: 
$$0\rightarrow\R\rightarrow C^{\infty}(M)\stackrel{{d}_{\Rc}^{\nabla^{E}}}{\rightarrow}\Rc^{1}(M;E)\stackrel{{d}_{\Rc}^{\nabla^{E}}}{\rightarrow}\ldots\stackrel{{d}_{\Rc}^{\nabla^{E}}}{\rightarrow}\Rc^{d}(M;E)\stackrel{{D}_{\Rc}^{\nabla^{E}}}{\rightarrow}\Rc^{d+1}(M;E)\stackrel{{d}_{\Rc}^{\nabla^{E}}}{\rightarrow}\ldots\stackrel{{d}_{\Rc}^{\nabla^{E}}}{\rightarrow}\Rc^{2d+1}(M;E)\rightarrow 0.$$

Note that ${{d}_{\Rc}^{\nabla^{E}}}^{2}\sgh=-{(R_{H}^{E}\wedge\sgh)}_{0}$ for $\sgh\in\Rc^{p}(M;E)\ (p\leq d-2)$. Also the previous sequence is not a complex excepted if $E$ is flat. Also we call this sequence the twisted Rumin pseudo-complex.\\

Let $N$ be a strictly pseudoconvex $CR$ manifold, $\phi:M\to N$ a horizontal map and $E=\ftn$ the pull-back bundle endowed with the connection $\nabla^{'}$ induced by the Tanaka-Wester connection of $TN$. Then $\dfh\in{\Omega}_{H}^{1}(M;\ftn)$ satisfies ${d}_{\Rc}^{\nabla^{'}}\dfh={({d}_{H}^{\nabla^{'}}\dfh)}_{0}=0$. Moreover, if $\phi$ is a pseudoharmonic map then $\delta_{\Rc}^{\nabla^{'}}\dfh=\delta_{H}^{\nabla^{'}}\dfh=0$. Consequently, if $d>1$, we have ${\triangle}_{\Rc}^{\nabla^{'}}\dfh=((d-1)\,{d}_{\Rc}^{\nabla^{'}}\delta_{\Rc}^{\nabla^{'}}+d\,\delta_{\Rc}^{\nabla^{'}}{d}_{\Rc}^{\nabla^{'}})\dfh=0$. If $d=1$ and $N$ is torsionless, then we have ${D}_{\Rc}^{\nabla^{'}}\dfh=\theta\wedge({\Lc}_{\xi}^{\nabla^{'}}\dfh+d_{H}^{\nabla^{'}}\delta_{H,J}^{\nabla^{'}}\dfh)=0$. The assumption $\phi$ pseudoharmonic yields to ${\triangle}_{\Rc}^{\nabla^{'}}\dfh={D}_{\Rc}^{\nabla^{'}*}{D}_{\Rc}^{\nabla^{'}}\dfh+{({d}_{\Rc}^{\nabla^{'}}{\delta}_{\Rc}^{\nabla^{'}})}^{2}\dfh=0$.  

\begin{Rem}
The condition $\phi:M\to N$ CR-pluriharmonic is equivalent to ${d}_{\Rc}^{\nabla^{'}}J^{*}\dfh=0$. Moreover, if $dim\,M=3$, it seems natural in view of the Theorem 5.1 to define the CR-pluriharmonicity of a horizontal map $\phi:M\to N$ by the condition ${D}_{\Rc}^{\nabla^{'}}J^{*}\dfh=0$. 
\end{Rem}

\section{Rigidity results for horizontal pseudoharmonic maps defined on contact locally sub-symmetric spaces}

Now we derive Mok-Siu-Yeung type formulas for horizontal maps from compact contact locally sub-symmetric spaces into strictly pseudoconvex $CR$ manifolds. In this section we assume that $M$ is a contact locally sub-symmetric space of dimension $2d+1\geq 5$. First we consider the case $M$ torsionless.

\begin{Le} Let $(M,\theta,\xi,J,g_{\theta})$ be a contact locally sub-symmetric space torsionless with $\sw$ non-zero. The tensor $Q_{H_{0}}^{+}\in\Gamma({S}^{2}({\wedge}_{H_{0}}^{2,+}(M)))$ given by $\displaystyle{Q_{H_{0}}^{+}=c_{0}I^{\C}_{H_{0}}+\Cm}$ with\\$\displaystyle{c_{0}=-\frac{8d}{d-1}\frac{\nn{\Cm}}{\sw}}$ is parallel and satisfies $\langle Q_{H_{0}}^{+},\rwo\rangle=0$ and ${(c_{H}(Q_{H_{0}}^{+}))}_{0}=0$. 
\end{Le}

\noindent{Proof.} First recall that $\displaystyle{I^{\C}_{H_{0}}=\frac{1}{4}{(g_{\theta_{H}}\kn g_{\theta_{H}})}_{0}^{+}}$. Now we determine $c_{0}$ in a such way that $\displaystyle{Q_{H_{0}}^{+}=c_{0}I^{\C}_{H_{0}}+\Cm}$ satisfies $\langle Q_{H_{0}}^{+},\rwo\rangle=0$. We have
$$\langle Q_{H_{0}}^{+},\rwo\rangle=\frac{c_{0}}{4}tr_{H}\widehat{\rwo}+\langle\Cm,\rwo\rangle.$$
Now, we have 
$$\rwo=\frac{\sw}{d(d+1)}I^{\C}_{H_{0}}+\frac{1}{2(d+2)}{\Bigl(\riwo\kn g_{\theta_{H}}-\rowo\kn\omega_{\theta}\Bigr)}_{0}+\Cm.$$
Since $\displaystyle{tr_{H}\widehat{{\Bigl(\riwo\kn g_{\theta_{H}}-\rowo\kn\omega_{\theta}\Bigr)}_{0}}=tr_{H}\widehat{\Cm}=0}$ and $\displaystyle{tr_{H}\,\widehat{I^{\C}_{H_{0}}}=\frac{d^2-1}{2}}$, we deduce that $\displaystyle{tr_{H}\widehat{\rwo}=\frac{d-1}{2d}\sw}$. 
Using $\langle\Cm,\rwo\rangle=\nn{\Cm}$, we obtain that
$$\langle Q_{H_{0}}^{+},\rwo\rangle=c_{0}\frac{d-1}{8d}\sw+\nn{\Cm}.$$
By taking $\displaystyle{c_{0}=-\frac{8d}{d-1}\frac{\nn{\Cm}}{\sw}}$ we obtain that $\langle Q_{H_{0}}^{+},\rwo\rangle=0$. Now we have 
$\displaystyle{c_{H}(Q_{H_{0}}^{+})=c_{0}\frac{d}{2}\Bigl(1-\frac{1}{d^2}\Bigr)g_{\theta_{H}}}$ and 
then ${(c_{H}(Q_{H_{0}}^{+}))}_{0}=0$. Since $M$ is a contact locally sub-symmetric space, then $\nabla\rw=0$ and $\sw$ constant yield to $\nabla\Cm=0$. Hence we have $\nn{\Cm}$ and $c_{0}$ constant. The parallelism of $Q_{H_{0}}^{+}=0$ directly follows. $\Box$

\begin{Prop} For any horizontal map $\phi$ from a compact contact locally sub-symmetric space $M$, holonomy irreducible and torsionless, to a Sasakian manifold $N$, we have: 
\begin{eqnarray}\label{h31}
&&\int_{M}\frac{c_{0}}{2d}\nn{{(\nbfh)}_{0}^{+}}+\langle\stackrel{\circ}{\Cm}{(\nbfh)}_{0}^{+},{(\nbfh)}_{0}^{+}\rangle
+\frac{c_{0}}{2}\Bigl(1-\frac{1}{d}\Bigr)\nn{{(\nbfh)}^{-}}\nonumber\\
&&+\langle\stackrel{\circ}{\Cm}{(\nbfh)}^{-},{(\nbfh)}^{-}\rangle-\frac{c_{0}}{2}\Bigl(1-\frac{1}{d^2}\Bigr)\nn{\delfh}-\frac{c_{0}}{2}(d^2-1)\nn{\dfxh}v_{g_{\theta}}\nonumber\\
&&=2\int_{M}c_{0}\Bigl(\frac{1}{d}\rf{2}{0}+\Bigl(1-\frac{1}{d}\Bigr)\rf{1}{1}\Bigr)+4\langle\Cm,\rwpf\rangle v_{g_{\theta}},
\end{eqnarray}
where $\displaystyle{c_{0}=-\frac{8d}{d-1}\frac{\nn{\Cm}}{\sw}}$. 
\end{Prop}

\noindent{Proof.} Let $\displaystyle{Q_{H_{0}}^{+}=c_{0}I^{\C}_{H_{0}}+\Cm}$ with $c_{0}$ defined in Lemma 6.1. The horizontal $J$-invariant symmetric $2$-tensor $(c_{H}(\widehat{\rw}\circ\widehat{Q_{H_{0}}^{+}}){)}^{\Sc}$ is parallel. We deduce from the irreducibility of $M$ that $(c_{H}(\widehat{\rw}\circ\widehat{Q_{H_{0}}^{+}}){)}^{\Sc}=\lambda g_{\theta_{H}}$ with $\lambda\in C^{\infty}(M,\R)$. Now,
$$\lambda=\frac{1}{2d}tr_{H}(c_{H}(\widehat{\rw}\circ\widehat{Q_{H_{0}}^{+}}){)}^{\Sc}=\frac{2}{d}tr_{H}(\widehat{\rw}\circ\widehat{Q_{H_{0}}^{+}})=\frac{4}{d}\langle\rw,Q_{H_{0}}^{+}\rangle=\frac{4}{d}\langle\rwo,Q_{H_{0}}^{+}\rangle=0.$$
Hence $(c_{H}(\widehat{\rw}\circ\widehat{Q_{H_{0}}^{+}}){)}^{\Sc}=0$. If $\phi$ is a horizontal map from $M$ to $N$, we have by Lemma 4.2\begin{eqnarray*}
\langle\stackrel{\circ}{Q_{H_{0}}^{+}}{(\nbfh)}_{0},{(\nbfh)}_{0}\rangle&=&\frac{c_{0}}{2}\Bigl(\frac{1}{d}\nn{{(\nbfh)}_{0}^{+}}+\Bigl(1-\frac{1}{d}\Bigr)\nn{{(\nbfh)}^{-}}\Bigr)\\
&+&\langle\stackrel{\circ}{\Cm}{(\nbfh)}_{0},{(\nbfh)}_{0}\rangle\\
tr_{H}\widehat{Q_{H_{0}}^{+}}&=&\frac{c_{0}}{2}\Bigl(d^2-1\Bigr)\\
\langle Q_{H_{0}}^{+},\rwpf\rangle&=&\frac{c_{0}}{4}\Bigl(\frac{1}{d}\rf{2}{0}+\Bigl(1-\frac{1}{d}\Bigr)\rf{1}{1}\Bigr)+\langle\Cm,\rwpf\rangle.
\end{eqnarray*}
By replacing in (\ref{h17}), we obtain the formula. $\Box$\\

Now we consider the case of contact locally sub-symmetric spaces with torsion. Let $M$ be a contact locally sub-symmetric spaces with torsion, we recall that we have $\displaystyle{{\tau}^{2}=\frac{|\tau|^{2}}{2d}id_{H}}$. We always may assume that $\displaystyle{\frac{|\tau|^{2}}{2d}=1}$ also $\tau$ becomes a paracomplex structure on $H$. Now $(\tau,J\circ\tau,J)$ defines a so called bi-paracomplex stucture (cf. \cite{ES}) on $M$. Any horizontal $2$-tensor $t_{H}$ on $M$ decomposes into $t_{H}=t_{H_{+}}+t_{H_{-}}$, where $t_{H_{\pm}}:=\frac{1}{2}(t_{H}\pm\tau^{*}t_{H})$ are respectively the $\tau$-invariant part and the $\tau$-anti-invariant part of $t_{H}$.
Let $\whtpm{\pm}$ be the bundle of $\tau$-(anti)invariant horizontal antisymmetric $2$-tensors.
For $Q_{H}\in S^{2}(\wh{2})$, we define $Q_{H_{\pm}}\in S^{2}(\whtpm{\pm})$ by 
\begin{eqnarray*} 
Q_{H_{\pm}}(X,Y,Z,W)&=&\frac{1}{4}(Q_{H}(X,Y,Z,W)\pm Q_{H}(\tau(X),\tau(Y),Z,W)\pm Q_{H}(X,Y,\tau(Z),\tau(W))\\
&+&Q_{H}(\tau(X),\tau(Y),\tau(Z),\tau(W))).
\end{eqnarray*}
The tensors ${(Q_{H})}^{\pm}_{\pm}\in S^{2}({(\wh{2})}^{\pm}_{\pm})$ are defined by ${(Q_{H})}^{\pm}_{\pm}={(Q_{H}^{\pm})}_{\pm}={(Q_{H_{\pm}})}^{\pm}.$

\begin{Le} Let $M$ be a contact locally sub-symmetric space with torsion and $s_{H}\in\sh{2}\otimes E$. We have the relations:
\begin{eqnarray*}
\stackrel{\circ}{A_{\theta}\kn A_{\theta}}s_{H}&=&2(\tau^{*}s_{H}-A_{\theta}\otimes tr_{H}c_{H}(A_{\theta}\otimes s_{H})),\\
\stackrel{\circ}{B_{\theta}\kn B_{\theta}}s_{H}&=&2({(J\circ\tau)}^{*}s_{H}-B_{\theta}\otimes tr_{H}c_{H}(B_{\theta}\otimes s_{H})),\\
c_{H}(A_{\theta}\kn A_{\theta})&=&c_{H}(B_{\theta}\kn B_{\theta})=-\frac{|\tau|^{2}}{d}g_{\theta_{H}},\quad
tr_{H}\,\widehat{A_{\theta}\kn A_{\theta}}=tr_{H}\,\widehat{B_{\theta}\kn B_{\theta}}=-|\tau|^{2}.
\end{eqnarray*}
\end{Le}
\newpage
\begin{Prop} For any horizontal map $\phi$ from a compact contact locally sub-symmetric space with torsion $M$ to a Sasakian manifold $N$, we have:
\begin{eqnarray*}
&&\int_{M}\nn{{(\nbfho)}^{+}_{+}}+\nn{{(\nbfh)}^{+}_{-}}+\nn{{(\nbfh)}^{-}_{+}+\frac{1}{d}A_{\theta}\otimes\delta^{\nabla^{'}}_{H}(\dfh\circ\tau)}\\
&&-\nn{{(\nbfh)}^{-}_{-}+\frac{1}{d}B_{\theta}\otimes\delta^{\nabla^{'}}_{H}(\dfh\circ J\circ\tau)}\\
&&+\Bigl(1-\frac{1}{d}\Bigr)\Bigl(\nn{\delta^{\nabla^{'}}_{H}(\dfh\circ J\circ\tau)}-\nn{\delta^{\nabla^{'}}_{H}(\dfh\circ\tau)}+d^{2}\nn{\dfxh}-\nn{\delfh}\Bigr)v_{g_{\theta}}\\
&&=8\int_{M}tr_{H}\widehat{{(\rwpf)}^{-}_{+}}v_{g_{\theta}},
\end{eqnarray*}
\begin{eqnarray}\label{h32}
&&\int_{M}\Bigl(1+\frac{2}{d}\Bigr)\nn{{(\nbfh)}^{+}_{-}}+\Bigl(1-\frac{2}{d}\Bigr)\nn{{(\nbfh)}^{-}_{+}+\frac{1}{d}A_{\theta}\otimes\delta^{\nabla^{'}}_{H}(\dfh\circ\tau)}\nonumber\\
&&+\Bigl(1-\frac{2}{d}\Bigr)\nn{{(\nbfh)}^{-}_{-}+\frac{1}{d}B_{\theta}\otimes\delta^{\nabla^{'}}_{H}(\dfh\circ J\circ\tau)}-\Bigl(1-\frac{2}{d}\Bigr)\nn{{(\nbfho)}^{+}_{+}}\nonumber\\
&&\Bigl(1-\frac{1}{d}\Bigr)\Bigl(1+\frac{2}{d}\Bigr)\Bigl(\nn{\delta^{\nabla^{'}}_{H}(\dfh\circ J\circ\tau)}+\nn{\delta^{\nabla^{'}}_{H}(\dfh\circ\tau)}-d^{2}\nn{\dfxh}-\nn{\delfh}\Bigr)v_{g_{\theta}}\nonumber\\
&&=8\int_{M}\Bigl(1-\frac{1}{d}\Bigr)tr_{H}\rwpfh{+}+\frac{1}{d}tr_{H}\rwpfh{-}-tr_{H}\widehat{{(\rwpf)}^{+}_{+}}v_{g_{\theta}},
\end{eqnarray}
\end{Prop}

\noindent{Proof.} Let the tensors $Q_{H}^{-}\in\Gamma({S}^{2}(\whpm{-})$ and $Q_{H_{0}}^{+}\in\Gamma({S}^{2}({\wedge}_{H_{0}}^{2,+}(M)))$ defined by:
$$Q_{H}^{-}={(g_{\theta_{H}}\kn g_{\theta_{H}})}^{-}_{+}=\frac{1}{4}(g_{\theta_{H}}\kn g_{\theta_{H}}-\omega_{\theta}\kn\omega_{\theta}+A_{\theta}\kn A_{\theta}-B_{\theta}\kn B_{\theta})$$
and
$$Q_{H_{0}}^{+}=\frac{1}{4}{({(g_{\theta_{H}}\kn g_{\theta_{H}})}_{0})}_{-}^{+}={(I^{\C}_{H_{0}})}_{-}=\frac{1}{2}(I^{\C}_{H_{0}}-\Tc_{H_{0}}),$$
with $\displaystyle{\Tc_{H_{0}}=\frac{1}{8}\Bigl(A_{\theta}\kn A_{\theta}+B_{\theta}\kn B_{\theta}+\frac{2}{d}\omega_{\theta}\odot\omega_{\theta}\Bigr)}$. Since $M$ is a contact locally sub-symmetric space, we have $\nabla_{H}A_{\theta}=\nabla_{H}B_{\theta}=0$ and then $\nabla_{H}Q_{H}^{-}=\nabla_{H}Q_{H_{0}}^{+}=0$. From Lemmas 4.2 and 6.2, we have 
$$c_{H}(Q_{H}^{-})=(d-1)g_{\theta_{H}},\quad c_{H}(Q_{H_{0}}^{+})=\frac{(d-1)(d+2)}{4d}g_{\theta_{H}},$$
and 
$$tr_{H}\widehat{Q_{H}^{-}}=d(d-1),\quad tr_{H}\widehat{Q_{H_{0}}^{+}}=\frac{(d-1)(d+2)}{4}.$$
It directly follows that ${(c_{H}(Q_{H}^{-}))}_{0}={(c_{H}(Q_{H_{0}}^{+}))}_{0}=0$, and that, 
$\displaystyle{\frac{tr_{H}\widehat{Q_{H}^{-}}}{d}B_{\theta}-\stackrel{\circ}{Q_{H}^{-}}B_{\theta}=0}$. Now since $\displaystyle{\rwo=\frac{2\sw}{d^2}{(I^{\C}_{H_{0}})}_{+}}$ then 
$$\widehat{\rw}\circ\widehat{Q_{H_{0}}^{+}}=\widehat{\rwo}\circ\widehat{Q_{H_{0}}^{+}}=\frac{2\sw}{d^2}\widehat{{(I^{\C}_{H_{0}})}_{+}}\circ\widehat{{(I^{\C}_{H_{0}})}_{-}}=0.$$
Also $(c_{H}(\widehat{\rwo}\circ\widehat{Q_{H_{0}}^{+}}){)}^{\Sc}=0$. Now, let $\phi$ be a horizontal map from $M$ to $N$, by Lemmas 4.2 and 6.2, we have 
\begin{eqnarray*}
\langle\stackrel{\circ}{Q_{H}^{-}}{(\nbfh)}_{0},{(\nbfh)}_{0}\rangle&=&\nn{{(\nbfho)}^{+}_{+}}+\nn{{(\nbfh)}^{+}_{-}}+\nn{{(\nbfh)}^{-}_{+}}-\nn{{(\nbfh)}^{-}_{-}}\\
&&-\frac{1}{4}\nn{tr_{H}c_{H}(A_{\theta}\otimes\nbfh)}+\frac{1}{4}\nn{tr_{H}c_{H}(B_{\theta}\otimes\nbfh)}\\
\langle\stackrel{\circ}{Q_{H_{0}}^{+}}{(\nbfh)}_{0},{(\nbfh)}_{0}\rangle&=&\frac{1}{4}\Bigl(\Bigl(1+\frac{2}{d}\Bigr)\nn{{(\nbfh)}^{+}_{-}}+\Bigl(1-\frac{2}{d}\Bigr)\nn{{(\nbfh)}^{-}_{+}}\\
&&+\Bigl(1-\frac{2}{d}\Bigr)\nn{{(\nbfh)}^{-}_{-}}-\Bigl(1-\frac{2}{d}\Bigr)\nn{{(\nbfho)}^{+}_{+}}\\
&&+\frac{1}{4}\nn{tr_{H}c_{H}(A_{\theta}\otimes\nbfh)}+\frac{1}{4}\nn{tr_{H}c_{H}(B_{\theta}\otimes\nbfh)}\Bigr).
\end{eqnarray*}
Now we have 
$$tr_{H}c_{H}(A_{\theta}\otimes\nbfh)=-2\delta^{\nabla^{'}}_{H}(\dfh\circ\tau)\quad{\rm and}\quad tr_{H}c_{H}(B_{\theta}\otimes\nbfh)=-2\delta^{\nabla^{'}}_{H}(\dfh\circ J\circ\tau).$$ 
Moreover 
$$\nn{{(\nbfh)}^{-}_{+}}=\nn{{(\nbfh)}^{-}_{+}+\frac{1}{d}A_{\theta}\otimes\delta^{\nabla^{'}}_{H}(\dfh\circ\tau)}+\frac{1}{d}\nn{\delta^{\nabla^{'}}_{H}(\dfh\circ\tau)}$$
and
$$\nn{{(\nbfh)}^{-}_{-}}=\nn{{(\nbfh)}^{-}_{-}+\frac{1}{d}B_{\theta}\otimes\delta^{\nabla^{'}}_{H}(\dfh\circ J\circ\tau)}+\frac{1}{d}\nn{\delta^{\nabla^{'}}_{H}(\dfh\circ J\circ\tau)}.$$
Then 
\begin{eqnarray*}
\langle\stackrel{\circ}{Q_{H}^{-}}{(\nbfh)}_{0},{(\nbfh)}_{0}\rangle&=&\nn{{(\nbfho)}^{+}_{+}}+\nn{{(\nbfh)}^{+}_{-}}
+\nn{{(\nbfh)}^{-}_{+}+\frac{1}{d}A_{\theta}\otimes\delta^{\nabla^{'}}_{H}(\dfh\circ\tau)}\\
&&-\nn{{(\nbfh)}^{-}_{-}+\frac{1}{d}B_{\theta}\otimes\delta^{\nabla^{'}}_{H}(\dfh\circ J\circ\tau)}\\
&&\Bigl(1-\frac{1}{d}\Bigr)\Bigl(\nn{\delta^{\nabla^{'}}_{H}(\dfh\circ J\circ\tau)}-\nn{\delta^{\nabla^{'}}_{H}(\dfh\circ\tau)}\Bigr)\\
\langle\stackrel{\circ}{Q_{H_{0}}^{+}}{(\nbfh)}_{0},{(\nbfh)}_{0}\rangle&=&\frac{1}{4}\Bigl(\Bigl(1+\frac{2}{d}\Bigr)\nn{{(\nbfh)}^{+}_{-}}-\Bigl(1-\frac{2}{d}\Bigr)\nn{{(\nbfho)}^{+}_{+}}\\
&&+\Bigl(1-\frac{2}{d}\Bigr)\nn{{(\nbfh)}^{-}_{+}+\frac{1}{d}A_{\theta}\otimes\delta^{\nabla^{'}}_{H}(\dfh\circ\tau)}\\
&&+\Bigl(1-\frac{2}{d}\Bigr)\nn{{(\nbfh)}^{-}_{-}+\frac{1}{d}B_{\theta}\otimes\delta^{\nabla^{'}}_{H}(\dfh\circ J\circ\tau)}\\
&&+\Bigl(1-\frac{1}{d}\Bigr)\Bigl(1+\frac{2}{d}\Bigr)\Bigl(\nn{\delta^{\nabla^{'}}_{H}(\dfh\circ J\circ\tau)}+\nn{\delta^{\nabla^{'}}_{H}(\dfh\circ\tau)}\Bigr)\Bigr).
\end{eqnarray*}
We have also
\begin{eqnarray*}
\langle Q_{H}^{-},\rwpf\rangle&=&\langle{(g_{\theta_{H}}\kn g_{\theta_{H}})}^{-}_{+},\rwpf\rangle=tr_{H}\widehat{{(\rwpf)}^{-}_{+}}\\
\langle Q_{H_{0}}^{+},\rwpf\rangle&=&\frac{1}{4}\langle{(g_{\theta_{H}}\kn g_{\theta_{H}})}^{+}_{-}-\frac{1}{d}\omega_{\theta}\odot\omega_{\theta},\rwpf\rangle\\
&=&\frac{1}{4}\Bigl(\Bigl(1-\frac{1}{d}\Bigr)tr_{H}\widehat{{(\rwpf)}^{+}_{-}}-\frac{1}{d}tr_{H}\widehat{{(\rwpf)}^{+}_{+}}+\frac{1}{d}tr_{H}\rwpfh{-}\Bigr)\\
&=&\frac{1}{4}\Bigl(\Bigl(1-\frac{1}{d}\Bigr)tr_{H}\rwpfh{+}+\frac{1}{d}tr_{H}\rwpfh{-}-tr_{H}\widehat{{(\rwpf)}^{+}_{+}}\Bigr).
\end{eqnarray*} 
By replacing in (\ref{h17}) and (\ref{h18}), we obtain the formulas. $\Box$\\

Now we deduce some rigidity results for the contact sub-symmetric space of non-compact type.\\

In the following we denote respectively by $\mathfrak{g},\mathfrak{k},\mathfrak{l}$ the Lie algebras of the Lie groups $G,K,L$. Let $\tilde{M}=G/K$ be a simply-connected contact sub-symmetric space of non-compact Hermitian type. Then $\tilde{M}$ is the total space of a $S^{1}$-fibration $\pi$ over an irreducible Hermitian symmetric space of non-compact type $\tilde{B}=G/L$. To the Hermitian symmetric space $G/L$, it is naturally associated an irreducible Hermitian orthogonal involutive Lie algebra $(\mathfrak{g},s,{\beta}_{/\mathfrak{p}})$, where $s$ is an involutive automorphism of $\mathfrak{g}$ such that $\mathfrak{l}$=$+1$-eigenspace of $s$, $\mathfrak{p}$=$-1$-eigenspace of $s$ and ${\beta}_{/\mathfrak{p}}$ is the $ad_{\mathfrak{l}}$-invariant inner product on $\mathfrak{p}$ given by the restriction of the Killing form of $\mathfrak{g}$ to $\mathfrak{p}$ (we refer to Falbel-Gorodski\cite{FG} for the precise definition). Also it follows a so-called irreducible subtorsionless Hermitian sub-orthogonal involutive Lie algebra $(\mathfrak{g},s,\mathfrak{k},{\beta}_{/\mathfrak{p}})$ associated to the sub-symmetric space $G/K$. Concerning $(\mathfrak{g},s,\mathfrak{k},{\beta}_{/\mathfrak{p}})$, the following facts hold \cite{FG}:\\
$\mathfrak{k}=[\mathfrak{l},\mathfrak{l}]$, $[\mathfrak{p},\mathfrak{p}]=\mathfrak{l}$ and $\mathfrak{l}=<\xi^{*}>\oplus\mathfrak{k}$ with $\xi^{*}$ in the center of $\mathfrak{l}$.\\
The ideal $\mathfrak{k}$ of $\mathfrak{l}$ is either a simple ideal or $\mathfrak{k}=\mathfrak{k_{1}}\oplus\mathfrak{k_{2}}$ with $\mathfrak{k_{1}}$ and $\mathfrak{k_{2}}$ are simple ideals of $\mathfrak{l}$.\\
The Killing form $\beta$ is negative definite on $\mathfrak{l}$ and we have the orthogonal decomposition of $\mathfrak{g}$ relatively to $\beta$, 
$$\mathfrak{g}=<\xi^{*}>\oplus\mathfrak{k}\oplus\mathfrak{p}.$$
The endomorphism $J^{*}={ad_{\xi^{*}}}_{/\mathfrak{p}}$ of $\mathfrak{p}$ defines a $ad_{\mathfrak{l}}$-invariant complex structure on $\mathfrak{p}$.\\
The $ad_{\mathfrak{l}}$-invariant skew-symmetric bilinear form $\beta(J^{*}., .)$ on $\mathfrak{p}$ is non-degenerate and coincides with $\Pi^{*}(d\theta)$ where $\Pi:G\to G/K$ is the natural projection.\\
Now the curvature ${\tilde{R}}_{H}^{\Wc}$ of $\tilde{M}$ is given by ${\tilde{R}}_{H}^{\Wc}=\pi^{*}R^{\tilde{B}}$ where $R^{\tilde{B}}$ is the curvature of $\tilde{B}$. Also $\tilde{M}$ has nonpositive pseudo-Hermitian sectional curvature and negative pseudo-Hermitian scalar curvature. The Lie algebra expression of ${\tilde{R}}_{H}^{\Wc}$ is given (cf. \cite{FGR}), for any $X_{i}^{*}\in\mathfrak{p}$, by: 
$${\tilde{R}}_{H}^{\Wc}(d\Pi(X_{1}^{*}),d\Pi(X_{2}^{*}),d\Pi(X_{3}^{*}),d\Pi(X_{4}^{*}))=\beta([X_{1}^{*},X_{2}^{*}],[X_{3}^{*},X_{4}^{*}]).$$Let $\displaystyle{c^{'}_{0}=-4\frac{\nn{{\tilde{R}}_{H}^{\Wc}}}{{\tilde{s}}^{\Wc}}=-4\frac{\nn{R^{\tilde{B}}}}{s^{\tilde{B}}}}(>0)$ and $\kappa(\tilde{M})$ be the lowest eigenvalue of the quadratic form $s_{H_{0}}\to\langle\stackrel{\circ}{{\tilde{R}}_{H}^{\Wc}}s_{H_{0}},s_{H_{0}}\rangle=\langle\stackrel{\circ}{\pi^{*}R^{\tilde{B}}}s_{H_{0}},s_{H_{0}}\rangle$ associated to ${\tilde{R}}_{H}^{\Wc}$ for horizontal traceless symmetric $2$-tensors. The following tabular, coming from those obtained in \cite{AB},\cite{CV} and \cite{YM} for the irreducible Hermitian symmetric spaces of non-compact type, gives the values of $c^{'}_{0}$ and $\kappa(\tilde{M})$ for the simply-connected contact sub-symmetric spaces of non-compact Hermitian type.\\

\begin{center}
\begin{tabular}{|c|c|c|c|}
\hline
Type & $d$ & $c^{'}_{0}$ & $\kappa(\tilde{M})$\\\hline
$SU(p,q)/SU(p)\times SU(q)$ & $pq$ & $\frac{pq+1}{{(p+q)}^{2}}$ & $-\frac{1}{p+q}$\\\hline
$SO^{*}(2p)/SU(p)$  & $\frac{p(p-1)}{2}$ & $\frac{1}{4}+\frac{3-p}{4{(p-1)}^{2}}$ & $-\frac{1}{2(p-1)}$\\\hline 
$Sp(p,\R)/SU(p)$ & $\frac{p(p+1)}{2}$ & $\frac{1}{4}+\frac{3+p}{4{(p+1)}^{2}}$ & $-\frac{1}{p+1}$\\\hline
$SO_{0}(p,2)/SO(p)$  & $p$ & $\frac{3}{2p}-\frac{1}{p^{2}}$ & $-\frac{1}{p}$\\\hline 
$E_{6(-14)}/Spin(10)$       & $16$ & $\frac{3}{16}$ & $-\frac{1}{12}$\\\hline 
$E_{7(-25)}/E_{6}$      & $27$ & $\frac{29}{162}$ & $-\frac{1}{18}$\\\hline 
\end{tabular}
\end{center}

 \begin{Th} Let $\tilde{M}$ be a simply-connected contact sub-symmetric space of non-compact Hermitian type other than $SU(d,1)/SU(d)$ and let $\Gamma$ be a cocompact discrete subgroup of $PsH(\tilde{M})$. Any horizontal pseudoharmonic map $\phi$ from $M=\tilde{M}/\Gamma$ to a Sasakian manifold $N$ with nonpositive pseudo-Hermitian complex sectional curvature satisfies $\nabla d\phi=0$.
\end{Th}

\noindent{Proof.} Let $\phi$ be a horizontal pseudoharmonic map from $M$ to $N$. Since $M$ and $N$ are torsionless, then (\ref{h25}) together with the assumption on the curvature of $N$ yields to $\dfxh=0$, ${(\nbfh)}_{0}^{+}=0$ and $\rf{2}{0}=0$. Now, using the irreducibility of $M$, we have by equation (\ref{h31}): 
\begin{eqnarray}\label{h33}
&&\int_{M}\frac{c_{0}}{2}\Bigl(1-\frac{1}{d}\Bigr)\nn{{(\nbfh)}^{-}}+\langle\stackrel{\circ}{\Cm}{(\nbfh)}^{-},{(\nbfh)}^{-}\rangle v_{g_{\theta}}\nonumber\\
&&=2\int_{M}c_{0}\Bigl(1-\frac{1}{d}\Bigr)\rf{1}{1}+4\langle\Cm,\rwpf\rangle v_{g_{\theta}}.
\end{eqnarray}
Now, since $M$ is pseudo-Einstein, then $\displaystyle{\Cm=\rw-\frac{\sw}{d(d+1)}I^{\C}_{H}}$ and $\displaystyle{\nn{\Cm}=\nn{\rw}-\frac{{\sw}^{2}}{4d(d+1)}}$. We deduce by Lemma 4.2 that
\begin{equation}\label{h34}
\frac{c_{0}}{2}\Bigl(1-\frac{1}{d}\Bigr)\nn{{(\nbfh)}^{-}}+\langle\stackrel{\circ}{\Cm}{(\nbfh)}^{-},{(\nbfh)}^{-}\rangle=c^{'}_{0}\nn{{(\nbfh)}^{-}}+\langle\stackrel{\circ}{\rw}{(\nbfh)}^{-},{(\nbfh)}^{-}\rangle
\end{equation}
and 
\begin{equation}\label{h35}
c_{0}\Bigl(1-\frac{1}{d}\Bigr)\rf{1}{1}+4\langle\Cm,\rwpf\rangle=2c^{'}_{0}\rf{1}{1}+4\langle\rw,\rwpf\rangle.
\end{equation}
By the comparison between $c^{'}_{0}$ and $\kappa(\tilde{M})$, we deduce that (\ref{h34}) is positive excepted if ${(\nbfh)}^{-}=0$. Moreover, using the Lie algebra expression of ${\tilde{R}}_{H}^{\Wc}$, we can prove as Jost-Yau do in \cite{JY} p 257-273, that (\ref{h35}) is always nonpositive. Also it follows from (\ref{h33}) that ${(\nbfh)}^{-}=0$. Consequently $\nbfh=0$ and $\dhfh=0$. The end of the proof follows from the proof of Theorem 5.2. $\Box$\\
 
Let $I$ be an open interval of $\R$ containing $0$, recall that a regular curve $c:I\to M$ on a strictly pseudoconvex $CR$ manifold $M$ is called a parabolic geodesic (cf. \cite{SD}) if $\dot{c}(0)\in H_{c(0)}$ and if there exists $\alpha\in\R$ such that
$\nabla_{\dot{c}(t)}\dot{c}(t)=\alpha\xi_{c(t)}$ for $t\in I$. As a consequence of the previous theorem, we have: 

\begin{Cor} Let $M=\tilde{M}/\Gamma$ as above. Any horizontal pseudoharmonic map $\phi$ from $M$ to a Sasakian manifold $N$ with nonpositive pseudo-Hermitian complex sectional curvature maps parabolic geodesics of $M$ to parabolic geodesics of $N$.
\end{Cor}

\noindent{Proof.} For any curve $c:I\subset\R\to M$ and any map $\phi$ from $M$ to $N$, we have 
$$\nabla_{\dot{(\phi\circ c)}(t)}^{'}\dot{(\phi\circ c)}(t)=(\nabla_{\dot{c}(t)}d\phi)(\dot{c}(t))+d\phi(\nabla_{\dot{c}(t)}\dot{c}(t)).$$If $\phi$ is a horizontal pseudoharmonic map from $M$ to $N$, we have $\nabla d\phi=0$. Consequently, for a parabolic geodesic $c:I\subset\R\to M$, we obtain that $\nabla_{\dot{(\phi\circ c)}(t)}^{'}\dot{(\phi\circ c)}(t)=\alpha d\phi(\xi_{c(t)})=\alpha f\xi_{(\phi\circ c)(t)}^{'}$. Hence $\phi\circ c$ is a parabolic geodesic of $N$. $\Box$

\begin{Cor} Let $M=\tilde{M}/\Gamma$ as above. Any horizontal pseudoharmonic map $\phi$ from $M$ to a Tanaka-Webster flat Sasakian manifold $N$ satisfies $\dfh=0$. 
\end{Cor}

\noindent{Proof.} Let $\phi$ be a horizontal pseudoharmonic map from $M$ to $N$. Since $N$ is Tanaka-Webster flat then Theorem 6.1 together with equation (\ref{h26}), yields to $\langle\dfh\circ\stackrel{\circ}{\riw},\dfh\rangle=0$. Now, since $M$ is pseudo-Einstein with $\sw<0$, then $\dfh=0$. $\Box$\\

If $M$ is a compact strictly pseudoconvex $CR$ manifold, then $b_{1}(M)=dim\,Ker\,{\triangle}_{\Rc}$.  Also, we have using formulas similar to (\ref{h25}),(\ref{h26}) and (\ref{h31}):
  
\begin{Cor} Let $\tilde{M}$ be a simply-connected contact sub-symmetric space of non-compact Hermitian type other than $SU(d,1)/SU(d)$ then $b_{1}(\tilde{M}/\Gamma)=0$.
\end{Cor}

\begin{Rem} We can observe that the previous corollary directly follows from the two following facts. First
$M$ is the total space of a $S^{1}$-fibration over a compact irreducible Hermitian locally symmetric space of non-compact type $B$ 
and then $b_{1}(M)=b_{1}(B)$. Second by the Matsushima Theorem \cite{YM} we have $b_{1}(B)=0$.
\end{Rem}

\section{Rigidity results for CR maps defined on contact locally sub-symmetric spaces}

In this section we suppose that $N$ is a strictly pseudoconvex $CR$ manifold and $\phi$ is a CR map from $M$ to $N$.

\begin{Prop}(i) For any CR map $\phi$ from a compact contact locally sub-symmetric space torsionless $M$ to $N$, we have:
\begin{eqnarray}\label{h36}
&&\int_{M}\frac{c_{0}}{2}\Bigl(1-\frac{1}{d}\Bigr)\nn{{(\nbfh)}^{-}}+\langle\stackrel{\circ}{\Cm}{(\nbfh)}^{-},{(\nbfh)}^{-}\rangle-c_{0}(d^2-1)\nn{\dfxh}v_{g_{\theta}}\nonumber\\
&&=2\int_{M}c_{0}\Bigl(1-\frac{1}{d}\Bigr){(HBK_{\phi}^{'\mathcal{W}})}_{H}+4\langle\Cm,\rwpf\rangle v_{g_{\theta}}.
\end{eqnarray}
(ii) For any CR map $\phi$ from a compact contact locally sub-symmetric space with torsion $M$ to $N$, we have:
\begin{eqnarray}\label{h37}
&&\int_{M}\Bigl(1-\frac{2}{d}\Bigr)\Bigr(\nn{{(\nbfh)}^{-}_{+}+\frac{1}{d}A_{\theta}\otimes\delta^{\nabla^{'}}_{H}(\dfh\circ\tau)}+\nn{{(\nbfh)}^{-}_{-}+\frac{1}{d}B_{\theta}\otimes J^{'}\delta^{\nabla^{'}}_{H}(\dfh\circ\tau)}\Bigr)\nonumber\\
&&+2\Bigl(1-\frac{1}{d}\Bigr)\Bigl(1+\frac{2}{d}\Bigr)\Bigl(\nn{\delta^{\nabla^{'}}_{H}(\dfh\circ\tau)}-d^{2}\nn{\dfxh}\Bigr)v_{g_{\theta}}=8\int_{M}\Bigl(1-\frac{1}{d}\Bigr){(HBK_{\phi}^{'\mathcal{W}})}_{H}-{(K_{\phi}^{'\mathcal{W}})}_{H}v_{g_{\theta}},\nonumber\\
\end{eqnarray}
with 
$${(K_{\phi}^{'\mathcal{W}})}_{H}=\sum_{i,j\leq d}\rwp(\dfh(\ep{i}),\dfh(\ep{j}),\dfh(\ep{i}),\dfh(\ep{j})),$$
and
$${(HBK_{\phi}^{'\mathcal{W}})}_{H}=\sum_{i,j\leq d}\rwp(\dfh(\ep{i}),J^{'}\dfh(\ep{i}),\dfh(\ep{j}),J^{'}\dfh(\ep{j})).$$
\end{Prop}

\noindent{Proof.} Let $\phi$ be a CR map from $M$ to $N$ then $J^{'}\circ\dfh=\dfh\circ J$. Also we have $\delfh=-d J^{'}\dfxh$ by (\ref{h27}) and $\delta^{\nabla^{'}}_{H}(\dfh\circ J\circ\tau)=J^{'}\delta^{\nabla^{'}}_{H}(\dfh\circ\tau)$. Moreover,
$${(\phi^{*}g_{\theta^{'}})}_{H}(X,Y)={(\phi^{*}\omega_{\theta^{'}})}_{H}(X,JY)=f\omega_{\theta}(X,JY)=f g_{\theta_{H}}(X,Y)$$
and 
$${(\phi^{*}B_{\theta^{'}})}_{H}(JX,JY)=B_{\theta^{'}}(J^{'}\dfh(X),J^{'}\dfh(Y))=-B_{\theta^{'}}(\dfh(X),\dfh(Y))=-{(\phi^{*}B_{\theta^{'}})}_{H}(X,Y).$$
Hence we have ${(\phi^{*}g_{\theta^{'}})}_{H}=f g_{\theta_{H}}$ and ${(\phi^{*}B_{\theta^{'}})}_{H}^{+}=0$. It follows that for each $Q_{H_{0}}^{+}$ defined in Propositions 6.1 and 6.2, we have  
$$\langle(c_{H}(\widehat{\rw}\circ\widehat{Q_{H_{0}}^{+}}){)}^{\Sc},{(\phi^{*}g_{\theta^{'}})}_{H}\rangle=\frac{1}{2}ftr_{H}(c_{H}(\widehat{\rw}\circ\widehat{Q_{H_{0}}^{+}}){)}^{\Sc}=4\langle\rwo,Q_{H_{0}}^{+}\rangle=0,$$
and
$$\langle\stackrel{\circ}{Q_{H_{0}}^{+}}{(\phi^{*}B_{\theta^{'}})}_{H},{(\phi^{*}g_{\theta^{'}})}_{H}\rangle=\frac{1}{2}ftr_{H}(\stackrel{\circ}{Q_{H_{0}}^{+}}{(\phi^{*}B_{\theta^{'}})}_{H}^{-})=0.$$
Now, we have ${(\phi^{*}\rwp)}_{H}^{-}=0$ and ${(\phi^{*}\rwp)}_{H}^{+}=\rwpf$. Consequently $tr_{H}\rwpfh{-}=0$ and $$tr_{H}\rwpfh{+}=\sum_{i,j\leq d}\rwp(\dfh(\ep{i}),J^{'}\dfh(\ep{i}),\dfh(\ep{j}),J^{'}\dfh(\ep{j})).$$  
We can choose an adapted frame $\{\ep{1},\ldots\ep{d},J\ep{1},\ldots J\ep{d}\}$ of $H$ such that 
$\tau(\ep{i})=\ep{i}$ and $\tau(J\ep{i})=-J\ep{i}$, also we obtain
\begin{eqnarray*}
tr_{H}\widehat{{(\rwpf)}^{+}_{+}}&=&\sum_{i,j\leq d}({(\rwpf)}_{+}(\ep{i},\ep{j},\ep{i},\ep{j})+{(\rwpf)}_{+}(\ep{i},J\ep{j},\ep{i},J\ep{j}))\\
&=&\sum_{i,j\leq d}\rwp(\dfh(\ep{i}),\dfh(\ep{j}),\dfh(\ep{i}),\dfh(\ep{j})).
\end{eqnarray*}
Now the formulas directly follow from equations (\ref{h17}) and (\ref{h18}) combinated to (\ref{h31}) and (\ref{h32}). $\Box$\\

A strictly pseudoconvex $CR$ manifold with constant holomorphic pseudo-Hermitian sectional curvature is called a pseudo-Hermitian space form. The sphere $S^{2d+1}=SU(d+1)/SU(d)$ viewed as the total space of the Hopf fibration over $\C P^{d}$ and its non-compact dual $SU(d,1)/SU(d)$, are examples of Sasakian pseudo-Hermitian space forms with respectively $\sw>0$ and $\sw<0$. The Heisenberg group $\Hc_{2d+1}$ with its standard pseudo-Hermitian structure is an example of flat Sasakian pseudo-Hermitian space form whereas the unit tangent bundle of the hyperbolic space $H^{d+1}$, $T_{1}H^{d+1}$ with its standard pseudo-Hermitian structure is an example of flat non-Sasakian pseudo-Hermitian space form (cf. \cite{JTC}). Note that all these examples are examples of contact sub-symmetric spaces. 

\begin{Th}(i) Let $\tilde{M}$ be a simply-connected contact sub-symmetric space of non-compact Hermitian type other than $SU(d,1)/SU(d)$ and let $\Gamma$ be a cocompact discrete subgroup of $PsH(\tilde{M})$. Then any horizontal pseudoharmonic CR map $\phi$ from $M=\tilde{M}/\Gamma$ to a strictly pseudoconvex $CR$ manifold $N$ with nonpositive pseudo-Hermitian complex sectional curvature satisfies $\nabla_{H}d\phi=0$.\\
(ii) Let $\tilde{M}$ be a simply-connected contact sub-symmetric space of non-compact type other than $\Hc_{2p+1}\times_{G}NCH$, then any horizontal pseudoharmonic CR map $\phi$ from $M=\tilde{M}/\Gamma$ to a pseudo-Hermitian space form $N$ with ${\sw}^{'}<0$ is constant.
\end{Th}

\noindent{Proof.} A horizontal pseudoharmonic CR map $\phi$ from $M$ to $N$ satisfies, by Proposition 5.2, $\delfh=\dfxh=0$ and ${(\nbfh)}^{+}=0$. Since $M$ is torsionless, then it directly follows from (\ref{h36}) that ${(\nbfh)}^{-}=0$. We deduce from the previous assumptions that $\nabla_{H}d\phi=0$. Hence (i) is proved. Now we assume that $N$ is a pseudo-Hermitian space form, then we have $\displaystyle{R^{\Wc}_{H^{'}}=\frac{{\sw}^{'}}{d^{'}(d^{'}+1)}I^{\C}_{H^{'}}}$. Since $\phi$ is a CR map, we obtain
$${(\phi^{*}\rwp)}_{H}^{+}=\rwpf=\frac{{\sw}^{'}}{d^{'}(d^{'}+1)}{(\phi^{*}I^{\C}_{H^{'}})}_{H}=f^{2}\frac{{\sw}^{'}}{d^{'}(d^{'}+1)}I^{\C}_{H}.$$
Hence we have $\langle\Cm,\rwpf\rangle=0$ and 
$${(HBK_{\phi}^{'\mathcal{W}})}_{H}=tr_{H}\rwpfh{+}=f^{2}\frac{{\sw}^{'}}{d^{'}(d^{'}+1)}tr_{H}\widehat{I^{\C}_{H}}=\frac{f^{2}}{2}\frac{d(d+1)}{d^{'}(d^{'}+1)}{\sw}^{'},$$
$${(K_{\phi}^{'\mathcal{W}})}_{H}=tr_{H}\widehat{{(\rwpf)}^{+}_{+}}=f^{2}\frac{{\sw}^{'}}{d^{'}(d^{'}+1)}tr_{H}\widehat{{(I^{\C}_{H})}_{+}}=\frac{f^{2}}{4}\frac{d(d-1)}{d^{'}(d^{'}+1)}{\sw}^{'},$$
and   
$$\Bigl(1-\frac{1}{d}\Bigr){(HBK_{\phi}^{'\mathcal{W}})}_{H}-{(K_{\phi}^{'\mathcal{W}})}_{H}=\frac{f^{2}}{4}\frac{(d-1)(d+2)}{d^{'}(d^{'}+1)}{\sw}^{'}.$$
Since ${\sw}^{'}<0$, the right hand sides of (\ref{h36}) and (\ref{h37}) are nonpositive, whereas the left hand sides are nonnegative.
We deduce from (\ref{h36}) that ${(HBK_{\phi}^{'\mathcal{W}})}_{H}=0$ and from (\ref{h37}) that $\displaystyle{\Bigl(1-\frac{1}{d}\Bigr){(HBK_{\phi}^{'\mathcal{W}})}_{H}-{(K_{\phi}^{'\mathcal{W}})}_{H}=0}$. In each case, we obtain that $f=0$. Since ${(\phi^{*}g_{\theta^{'}})}_{H}=f g_{\theta_{H}}$ and $d\phi(\xi)=f\xi^{'}$, then $\phi$ is constant. $\Box$\\

A regular curve $c:I\to M$ on a strictly pseudoconvex $CR$ manifold $M$ is called a Carnot-Caratheodory geodesic (cf. \cite{BD2},\cite{MR},\cite{RS}) if $\dot{c}(t)\in H_{c(t)}$ and if there exists a function $\alpha:I\to\R$ with $\dot{\alpha}(t)=A_{\theta}(\dot{c}(t),\dot{c}(t))$ such that
$\nabla_{{(\dot{c}(t))}_{H}}{(\dot{c}(t))}_{H}=-\alpha(t)J{(\dot{c}(t))}_{H}$ for $t\in I$. Also, we have 

\begin{Cor} Let $\tilde{M}$ be a simply-connected contact sub-symmetric space of non-compact Hermitian type other than $SU(d,1)/SU(d)$ and $M=\tilde{M}/\Gamma$. Then any horizontal pseudoharmonic CR map $\phi$ from $M$ to a Sasakian manifold $N$ with nonpositive pseudo-Hermitian complex sectional curvature maps Carnot-Caratheodory geodesics of $M$ to Carnot-Caratheodory geodesics of $N$.
\end{Cor}

\noindent{Proof.} Let $c:I\subset\R\to M$ be a Carnot-Caratheodory geodesic, then $\nabla_{{(\dot{c}(t))}_{H}}{(\dot{c}(t))}_{H}=-\alpha J{(\dot{c}(t))}_{H}$ with $\alpha$ constant (since $M$ is torsionless). Let $\phi$ be a horizontal pseudoharmonic CR map from $M$ to $N$, we obtain using $\nabla_{H} d\phi=0$ that $\nabla_{{(\dot{(\phi\circ c)}(t))}_{H^{'}}}({\dot{(\phi\circ c)}(t))}_{H^{'}}=-\alpha\dfh(J{(\dot{c}(t))}_{H})=-\alpha J^{'}\dfh({(\dot{c}(t))}_{H})=-\alpha J^{'}{(\dot{(\phi\circ c)}(t))}_{H^{'}}$. Also $\phi\circ c$ is a Carnot-Caratheodory geodesic of $N$. $\Box$

\end{document}